# A FUNCTIONAL LIMIT THEOREM FOR THE PROFILE OF SEARCH TREES

By Michael Drmota,[1] Svante Janson and Ralph Neininger[2]

*TU Wien, Uppsala University and J. W. Goethe University*

We study the profile $X_{n,k}$ of random search trees including binary search trees and $m$-ary search trees. Our main result is a functional limit theorem of the normalized profile $X_{n,k}/\mathbb{E}X_{n,k}$ for $k = \lfloor \alpha \log n \rfloor$ in a certain range of $\alpha$.

A central feature of the proof is the use of the contraction method to prove convergence in distribution of certain random analytic functions in a complex domain. This is based on a general theorem concerning the contraction method for random variables in an infinite-dimensional Hilbert space. As part of the proof, we show that the Zolotarev metric is complete for a Hilbert space.

**1. Introduction.** Search trees are used in computer science as data structures that hold data (also called keys) from a totally ordered set in order to support operations on the data such as searching and sorting. After having constructed the search tree for a set of keys, the complexity of operations performed on the data is identified by corresponding shape parameters of the search tree (examples are given below). Usually, one assumes a probabilistic model for the set of data or uses randomized procedures to build up search trees so that the resulting trees become random and the typical complexity of operations can be captured by computing expectations, variances, limit laws or tail bounds. In this paper, we study the profile of a general class of random search trees that includes many trees used in computer science such as the binary search tree and $m$-ary search trees with respect to functional limit laws.

A random binary search tree is constructed for a set of keys as follows. One key, the so-called *pivot*, is chosen uniformly from the set of data and

Received September 2006; revised June 2007.

[1]Supported by the Austrian Science Foundation FWF, project S9604.

[2]Supported by an Emmy Noether fellowship of the DFG.

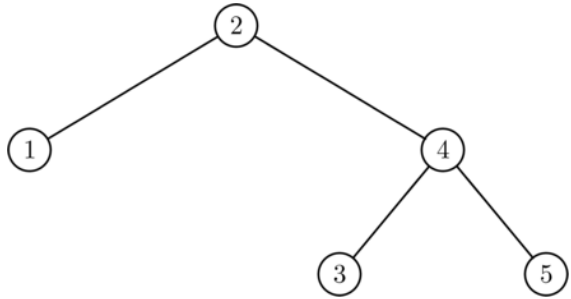

Fig. 1. *A random binary search tree for the data set* $\{1,2,3,4,5\}$. *In the first step, key 2 is chosen as the pivot. For the right subtree of the root holding the keys* $\{3,4,5\}$, *key 4 is chosen as a pivot. The profile of this tree is* $(1,2,2,0,0,\ldots)$.

inserted in the root of the tree. All other keys are compared with the pivot. Those which are smaller are used to build a random binary search tree as the left subtree of the root; those which are larger (or equal) than the pivot are used to build the right subtree of the root. For building these subtrees, the procedure is recursively applied. An example is given in Figure 1.

For the general class of search trees, explained in Section 2 and studied in this paper, this construction rule is generalized so that nodes may hold $m-1 \geq 1$ keys and have $m$ subtrees and, further, the rule to choose the pivots may be more general, resulting in more balanced trees as a parameter $t \geq 0$ is increased; see Section 2. For example, if $m=2$, then the pivot is chosen as the median of $2t+1$ random elements. This more general search tree model reduces to the binary search tree for the choice $(m,t)=(2,0)$.

The *depth* of a key in the tree is its node's distance to the root of the tree. This quantity is a measure of the complexity involved in searching for the number inserted in that node. Other quantities, important in the context of computer science, are the *internal path length* of the tree, which is the sum of the depths of all keys, and the *height* of the tree, which is the maximal depth in the tree.

In this paper, we study the profile of search trees, which is the infinite vector $\mathbf{X}_n = (X_{n,k})_{k \geq 0}$, where $X_{n,k}$ is the number of keys that are stored in nodes with depth $k$.

The profile of binary search trees (and related structures) has been intensively studied in the literature [4, 6, 7, 8, 10, 11, 12, 15, 17, 25]. Most results concern 1st and 2nd moments. However, there are also distributional results, particularly for binary search trees and recursive trees [4, 6, 15] that are of the form

$$\frac{X_{n,\lfloor \alpha \log n \rfloor}}{\mathbb{E} X_{n,\lfloor \alpha \log n \rfloor}} \xrightarrow{\mathrm{d}} X(\alpha)$$

for fixed $\alpha$ (contained in a suitable interval). The advantage of binary search trees and recursive trees is that there is an underlying martingale structure which also allows functional limit theorems to be proven (see [4, 6] for binary search trees). Unfortunately, this martingale structure is (generally) missing in the kind of trees that we want to study.

Our main result is the following, where we actually prove functional convergence of random functions on an interval $I'$. More precisely, we use the space $D(I')$ of right-continuous functions with left-hand limits equipped with the Skorohod topology; see Section 4 for the definition and note that when, as here, the limit is continuous, convergence in the Skorohod topology is equivalent to uniform convergence on every compact subinterval.

In the formulation of Theorem 1.1, we also use the function $\lambda_1(z)$, defined in Section 3 as the dominant root of (3.4), and the stochastic process $(Y(z), z \in B)$ (of analytic functions in a certain domain $B$ containing the interval $I$) that is defined as the unique solution of a stochastic fixed point equation (3.7) which is discussed in Section 9, satisfying the further conditions that $\mathbb{E}Y(z) = 1$ and that for each $x \in I$, there exists an $s(x) > 1$ such that $\mathbb{E}|Y(z)|^{s(x)}$ is finite and bounded in a neighborhood of $x$.

THEOREM 1.1.  *Let $m \geq 2$ and $t \geq 0$ be given integers and let $(X_{n,k})_{k \geq 0}$ be the profile of the corresponding random search tree with $n$ keys.*

*Set $I = \{\beta > 0 : 1 < \lambda_1(\beta^2) < 2\lambda_1(\beta) - 1\}$, $I' = \{\beta\lambda_1'(\beta) : \beta \in I\}$ and let $\beta(\alpha) > 0$ be defined by $\beta(\alpha)\lambda_1'(\beta(\alpha)) = \alpha$. We then have, in $D(I')$, that*

$$(1.1) \qquad \left( \frac{X_{n, \lfloor \alpha \log n \rfloor}}{\mathbb{E}X_{n, \lfloor \alpha \log n \rfloor}}, \alpha \in I' \right) \overset{\mathrm{d}}{\longrightarrow} (Y(\beta(\alpha)), \alpha \in I').$$

REMARK 1.1.  From the definitions of $I$ and $I'$, it is not clear that they are in fact intervals. We will make this precise in Lemma 8.5.

REMARK 1.2.  In exactly the same way, one can consider other similarly defined parameters. For example, in Section 11, we discuss the external profile.

The proof of Theorem 1.1 is divided into several steps. After defining suitable function spaces (Section 4), we show (Section 9) the following theorem, which states that if $W_n(z) := \sum_k X_{n,k} z^k$ are the profile polynomials, then the normalized profile polynomials $W_n(z)/\mathbb{E}W_n(z)$ converge weakly to $Y(z)$ for $z$ contained in a suitable complex region $B$, where $Y(z)$ is, as above, the solution of a stochastic fixed point equation (3.7). Note that convergence in $\mathcal{H}(B)$ means uniform convergence on every compact subset of $B$.

THEOREM 1.2. *There exists a complex region $B$ that contains the real interval $(1/m, \beta(\alpha_+))$, where $\alpha_+$ is defined in (1.3), and an analytic stochastic process $(Y(z), z \in B)$ satisfying (3.7) and $\mathbb{E}Y(z) = 1$, such that, in $\mathcal{H}(B)$,*

$$(1.2) \qquad \left( \frac{W_n(z)}{\mathbb{E}W_n(z)}, z \in B \right) \xrightarrow{\mathrm{d}} (Y(z), z \in B).$$

Finally, we apply a suitable continuous functional (which is related to Cauchy's formula) in order to derive Theorem 1.1 from this property (Section 10).

Important tools in this argument are Theorems 5.1 and 6.1, which show that one can use the contraction method with the Zolotarev metric $\zeta_s$ for random variables with values in a separable Hilbert space. (We do not know whether these theorems extend to arbitrary Banach spaces.)

In the special case of binary search trees, Theorems 1.1 and 1.2 have been proven earlier, also in stronger versions [4, 6, 7].

Before we go into the details, we wish to comment on the interval $I$ of Theorem 1.1. It is well known that the height of random search trees is of order $\log n$. Thus, it is natural that there might be a restriction on the parameter $\alpha = k/\log n$, where $k$ denotes the depth.

In fact, there are several *critical values* for $\alpha = k/\log n$, namely

- $\alpha = \alpha_0 := \left( \dfrac{1}{t+1} + \dfrac{1}{t+2} + \cdots + \dfrac{1}{(t+1)m - 1} \right)^{-1}$;

- $\alpha = \alpha_{\max} := \left( \dfrac{1}{t+2} + \dfrac{1}{t+3} + \cdots + \dfrac{1}{(t+1)m} \right)^{-1}$;

- $\alpha = \alpha_+$, where $\alpha_+ > \alpha_0$ is the solution of the equation

$$(1.3) \qquad \lambda_1(\beta(\alpha)) - \alpha \log(\beta(\alpha)) - 1 = 0.$$

In order to explain these critical values, we must look at the expected profile $\mathbb{E}X_{n,k}$. If $\alpha = k/\log n \le \alpha_0 - \varepsilon$ (for some $\varepsilon > 0$), then

$$\mathbb{E}X_{n,k} \sim (m-1)m^k,$$

whereas if $\alpha = k/\log n \ge \alpha_0 + \varepsilon$, then

$$\mathbb{E}X_{n,k} \sim \frac{E(\beta(\alpha))n^{\lambda_1(\beta(\alpha)) - \alpha \log(\beta(\alpha)) - 1}}{\sqrt{2\pi(\alpha + \beta(\alpha)^2 \lambda_1''(\beta(\alpha)))\log n}}$$

for some continuous function $E(z)$; see Lemma 8.3. This means that up to level $k = \alpha_0 \log n$, the tree is (almost) complete. Note that the critical value $k/\log n = \alpha_0$ corresponds to $z = \beta = 1/m$ and $\lambda_1(1/m) = 1$, and thus that

$$n^{\lambda_1(\beta(\alpha_0)) - \alpha_0 \log(\beta(\alpha_0)) - 1} = n^{\alpha_0 \log m} = m^k.$$

We can be even more precise. If $\alpha = k/\log n \in [\varepsilon, \alpha_0 - \varepsilon]$, then

$$\mathbb{E}X_{n,k} = (m-1)m^k - r_{n,k},$$

with

$$r_{n,k} \sim \frac{E_1(\beta(\alpha))n^{\lambda_1(\beta(\alpha)) - \alpha\log(\beta(\alpha)) - 1}}{\sqrt{2\pi(\alpha + \beta(\alpha)^2\lambda_1''(\beta(\alpha)))\log n}}$$

for some continuous function $E_1(z)$.

The second critical value $k/\log n = \alpha_{\max}$ corresponds to $z = \beta = 1$ and $\lambda_1(1) = 2$. Here, we have

$$\mathbb{E}X_{n,k} \sim \frac{n}{\sqrt{2\pi(\alpha_{\max} + \lambda_1''(1))\log n}} \exp\left(-\frac{(k - \alpha_{\max}\log n)^2}{2(\alpha_{\max} + \lambda_1''(1))\log n}\right)$$

[uniformly for $k = \alpha_{\max}\log n + O(\sqrt{\log n})$]. This means that most nodes are concentrated around that level. In fact, $\alpha_{\max}\log n$ is the expected depth.

Finally, if $\alpha = k/\log n < \alpha_+$, then $\mathbb{E}X_{n,k} \to \infty$ and if $\alpha = k/\log n > \alpha_+$, then $\mathbb{E}X_{n,k} \to 0$. This means that the range $\alpha = k/\log n \in (0, \alpha_+)$ is exactly the range where the profile $X_{n,k}$ is actually present.

We also see that the interval $I'$ of Theorem 1.1 is strictly contained in $(\alpha_0, \alpha_+)$, but we have $\alpha_{\max} \in I'$. This means that we definitely cover the *most important range*. However, it seems that Theorem 1.1 is not optimal. The condition $\lambda_1(\beta^2) < 2\lambda_1(\beta) - 1$ comes from the fact that we are using $L^2$ techniques in order to derive Theorem 1.1 from Theorem 1.2. We conjecture that this is just a technical restriction and that Theorem 1.1 actually holds for $\alpha \in (\alpha_0, \alpha_+)$.

Incidentally, $r_{n,k}$ has a similar critical value $\alpha_- < \alpha_0$ that is the second positive solution of (1.3). If $\alpha < \alpha_-$, then $r_{n,k} \to 0$ and if $\alpha > \alpha_-$, then $r_{n,k} \to \infty$. The two constants $\alpha_-, \alpha_+$ are related to the *speed* of the leftmost and rightmost particles in suitable discrete branching random walks (see [5]). Note that they can be also computed by

$$\alpha_- = \left(\sum_{j=0}^{(t+1)(m-1)-1} \frac{1}{\lambda_- + t + j}\right)^{-1}$$

and

$$\alpha_+ = \left(\sum_{j=0}^{(t+1)(m-1)-1} \frac{1}{\lambda_+ + t + j}\right)^{-1},$$

where $\lambda_-$ and $\lambda_+$ are the two solutions of

$$\sum_{j=0}^{(t+1)(m-1)-1} \log(\lambda + t + j) - \log(m(tm + m - 1)!/t!)$$

(1.4)

$$= \sum_{j=0}^{(t+1)(m-1)-1} \frac{\lambda-1}{\lambda+t+j}.$$

Further, the expected height of $m$-ary search trees satisfies $\mathbb{E}H_n \sim \alpha_+ \log n$ and the expected saturation level $\mathbb{E}\tilde{H}_n \sim \alpha_- \log n$.

NOTATION.  If $f$ and $g$ are two functions on the same domain, then $f \lesssim g$ means the same as $f = O(g)$, that is, $|f| \le Cg$ for some constant $C$.

**2. Random search trees.**  To describe the construction of the search tree, we begin with the simplest case $t = 0$. If $n = 0$, the tree is empty. If $1 \le n \le m-1$, the tree consists of a root only, with all keys stored in the root. If $n \ge m$, we randomly select $m-1$ keys that are called *pivots* (with the uniform distribution over all sets of $m-1$ keys). The pivots are stored in the root. The $m-1$ pivots split the set of the remaining $n-m+1$ keys into $m$ subsets $I_1, \ldots, I_m$: if the pivots are $x_1 < x_2 < \cdots < x_{m-1}$, then $I_1 := \{x_i : x_i < x_1\}$, $I_2 := \{x_i : x_1 < x_i < x_2\}, \ldots, I_m := \{x_i : x_{m-1} < x_i\}$. We then recursively construct a search tree for each of the sets $I_i$ of keys (ignoring $I_i$ if it is empty) and attach the roots of these trees as children of the root in the search tree.

In the case $m = 2$, $t = 0$, we thus have the well-studied *binary search tree* [4, 6, 7, 11, 12, 15, 26].

In the case $t \ge 1$, the only difference is that the pivots are selected in a different way, which affects the probability distribution of the set of pivots and thus of the trees. We now select $mt + m - 1$ keys at random, order them as $y_1 < \cdots < y_{mt+m-1}$ and let the pivots be $y_{t+1}, y_{2(t+1)}, \ldots, y_{(m-1)(t+1)}$. In the case $m \le n < mt + m - 1$, when this procedure is impossible, we select the pivots by some supplementary rule (possibly random, but depending only on the order properties of the keys); our results do not depend on the choice of this supplementary rule.

This splitting procedure was first introduced by Hennequin for the study of variants of the Quicksort algorithm and is referred to as the *generalized Hennequin Quicksort* (cf. Chern, Hwang and Tsai [9]).

In particular, in the case $m = 2$, we let the pivot be the median of $2t + 1$ randomly selected keys (when $n \ge 2t + 1$).

We describe the splitting of the keys by the random vector $\mathbf{V}_n = (V_{n,1}, V_{n,2}, \ldots, V_{n,m})$, where $V_{n,k} := |I_k|$ is the number of keys in the $k$th subset and thus the number of nodes in the $k$th subtree of the root (including empty subtrees).

We thus always have, provided $n \ge m$,

$$V_{n,1} + V_{n,2} + \cdots + V_{n,m} = n - (m-1) = n + 1 - m$$

and elementary combinatorics, counting the number of possible choices of the $mt + m - 1$ selected keys, shows that the probability distribution is, for $n \geq mt + m - 1$ and $n_1 + n_2 + \cdots + n_m = n - m + 1$,

$$(2.1) \qquad \mathbb{P}\{\mathbf{V}_n = (n_1, \ldots, n_m)\} = \frac{\binom{n_1}{t} \cdots \binom{n_m}{t}}{\binom{n}{mt+m-1}}.$$

(The distribution of $\mathbf{V}_n$ for $m \leq n < mt + m - 1$ is not specified.)

In particular, for $n \geq mt + m - 1$, the components $V_{n,j}$ are identically distributed and another simple counting argument yields, for $n \geq mt + m - 1$ and $0 \leq \ell \leq n - 1$,

$$(2.2) \qquad \mathbb{P}\{V_{n,j} = \ell\} = \frac{\binom{\ell}{t}\binom{n-\ell-1}{(m-1)t+m-2}}{\binom{n}{mt+m-1}}.$$

For example, for the binary search tree with $m = 2$ and $t = 0$, we thus have $V_{n,1}$ and $V_{n,2} = n - 1 - V_{n-1}$ uniformly distributed on $\{0, \ldots, n-1\}$.

**3. The profile polynomial.** The recursive construction of the random search tree in Section 2 leads to a recursion for the profile $\mathbf{X}_n = (X_{n,k})_{k \geq 0}$:

$$(3.1) \qquad X_{n,k} \stackrel{\mathrm{d}}{=} X^{(1)}_{V_{n,1},k-1} + X^{(2)}_{V_{n,2},k-1} + \cdots + X^{(m)}_{V_{n,m},k-1},$$

jointly in $k \geq 0$ for every $n \geq m$, where the random vector $\mathbf{V}_n = (V_{n,1}, V_{n,2}, \ldots, V_{n,m})$ is as in Section 2 and is the same for every $k \geq 0$, and $\mathbf{X}_n^{(j)} = (X_{n,k}^{(j)})_{k \geq 0}$, $j = 1, \ldots, m$, are independent copies of $\mathbf{X}_n$ that are also independent of $\mathbf{V}_n$. We further have $X_{n,0} = m - 1$ for $n \geq m$. For $n \leq m - 1$, we simply have $X_{n,0} = n$ and $X_{n,k} = 0$, $k \geq 1$.

Note that, by induction, $X_{n,k} = 0$ when $k \geq n$. Hence, each vector $\mathbf{X}_n$ has only a finite number of nonzero components.

Let $W_n(z) = \sum_k X_{n,k} z^k$ denote the random *profile polynomial*. By (3.1), it is recursively given by $W_n(z) = n$ for $n \leq m - 1$ and

$$(3.2) \quad W_n(z) \stackrel{\mathrm{d}}{=} zW^{(1)}_{V_{n,1}}(z) + zW^{(2)}_{V_{n,2}}(z) + \cdots + zW^{(m)}_{V_{n,m}}(z) + m - 1, \qquad n \geq m,$$

where $W^{(j)}_\ell(z)$, $j = 1, \ldots, m$, are independent copies of $W_\ell(z)$ that are independent of $\mathbf{V}_n$, $\ell \geq 0$. From this relation, we obtain a recurrence for the expected profile polynomial $\mathbb{E}W_n(z)$. We have, using (2.2), for $n \geq mt + m - 1$,

$$(3.3) \qquad \mathbb{E}W_n(z) = mz \sum_{\ell=0}^{n-1} \frac{\binom{\ell}{t}\binom{n-\ell-1}{(m-1)t+m-2}}{\binom{n}{mt+m-1}} \mathbb{E}W_\ell(z) + m - 1.$$

For any fixed complex $z$, this is a recursion of the type studied in Chern, Hwang and Tsai [9]. More precisely, it fits ([9], (13)) with $a_n = \mathbb{E}W_n(z)$, $r = mt + m - 1$ and $c_t = mzr!/t!$, while $c_j = 0$ for $j \neq t$. Further, $b_n = m - 1$ for $n \geq mt + m - 1$, while $b_n = a_n = \mathbb{E}W_n(z)$ for $n < mt + m - 1$.

It follows from [9] that the asymptotics of $\mathbb{E}W_n(z)$ as $n \to \infty$ depend on the roots of the indicial polynomial

$$
\begin{aligned}
\text{(3.4)} \quad \Lambda(\theta; z) &:= \theta^{\overline{mt+m-1}} - mz\frac{(mt+m-1)!}{t!}\theta^{\overline{t}} \\
&= \theta(\theta+1)\cdots(\theta+mt+m-2) \\
&\quad - mz\frac{(mt+m-1)!}{t!}\theta(\theta+1)\cdots(\theta+t-1)
\end{aligned}
$$

using the notation $x^{\overline{m}} := x(x+1)\cdots(x+m-1) = \Gamma(x+m)/\Gamma(x)$. If we set

$$
\text{(3.5)} \quad F(\theta) := \frac{t!}{m(mt+m-1)!}(\theta+t)(\theta+t+1)\cdots(\theta+mt+m-2),
$$

then

$$
\Lambda(\theta; z) = \frac{m(mt+m-1)!}{t!}\theta^{\overline{t}}(F(\theta) - z),
$$

which implies that the roots of $\Lambda(\lambda; z) = 0$ are $0, -1, -2, \ldots, -t+1$ (if $t \geq 1$) together with the roots of $F(\theta) = z$. Let $\lambda_j(z)$, $j = 1, \ldots, (m-1)(t+1)$, denote the roots of $F(\theta) = z$ (counted with multiplicities), arranged in decreasing order of their real parts: $\Re\lambda_1(z) \geq \Re\lambda_2(z) \geq \cdots$.

Further, let $D_s$, for real $s$, be the set of all complex $z$ such that $\Re\lambda_1(z) > s$ and $\Re\lambda_1(z) > \Re\lambda_2(z)$ [in particular, $\lambda_1(z)$ is a simple root]. It is easily seen that the set $D_s$ is open and that $\lambda_1(z)$ is an analytic function of $z \in D_s$. If $z \in D_s$ is real, then $\lambda_1(z)$ must be real (and thus greater than $s$) because otherwise, $\overline{\lambda_1(z)}$ would be another root with the same real part.

By [9], Theorem 1(i), we have the following result. Note that $K_0$ and $K_1$ [our $E(z)$] in [9], Theorem 1(i), are analytic functions of $z$ and $\lambda_1$, and thus of $z \in D_1$, and that they are positive for $\lambda_1 > 0$ because $b_k = m - 1 > 0$ for $k \geq mt + m - 1$ and $b_k = \mathbb{E}W_k(z) \geq 0$ for smaller $k$. (See also Lemma 8.2 and the Appendix.)

LEMMA 3.1.   *If $z \in D_1$, then*

$$
\mathbb{E}W_n(z) = (E(z) + o(1))n^{\lambda_1(z)-1}
$$

*for some analytic function $E(z)$ with $E(z) > 0$ for $z \in D_1 \cap (0, \infty)$.*

LEMMA 3.2.   *The set $D_1$ is an open domain in the complex plane that contains the interval $(1/m, \infty)$.*

[Lemma 3.2 will be proven in a more general context in Lemma 8.1. Note that $F(1) = 1/m$ and thus $\lambda_1(1/m) = 1$.]

Set $M_n(z) = W_n(z)/G_n(z)$, where $G_n(z) = \mathbb{E}W_n(z)$. Then (3.2) can be rewritten as

$$M_n(z) \stackrel{\mathrm{d}}{=} \frac{G_{V_{n,1}}(z)}{G_n(z)} z M_{V_{n,1}}^{(1)}(z) + \cdots + \frac{G_{V_{n,m}}(z)}{G_n(z)} z M_{V_{n,m}}^{(m)}(z) + \frac{m-1}{G_n(z)}.$$

Note that $G_V(z)$, where $V$ is an integer-valued random variable, is considered as the random variable $\mathbb{E}W_n(z)|_{n=V}$ and not as $\mathbb{E}W_V(z)$, that is, the expected value is only taken with respect to $\mathbf{X}_n$. Next, let the random vector $\mathbf{V} = (V_1, V_2, \ldots, V_m)$ be supported on the simplex $\Delta = \{(s_1, \ldots, s_m) : s_j \geq 0, s_1 + \cdots + s_m = 1\}$ with density

$$f(s_1, \ldots, s_m) = \frac{((t+1)m-1)!}{(t!)^m}(s_1 \cdots s_m)^t,$$

where $t \geq 0$ is the same integer parameter as above. (This is known as a *Dirichlet distribution.*) It is easy to show that

$$(3.6) \qquad \frac{1}{n}\mathbf{V}_n \stackrel{\mathrm{d}}{\longrightarrow} \mathbf{V} \qquad \text{as } n \to \infty.$$

REMARK 3.1. For $n \geq mt + m - 1$, the shifted random vector $(V_{n,1} - t, \ldots, V_{n,m} - t)$ has a multivariate Pólya–Eggenberger distribution that can be defined as the distribution of the vector of the numbers of balls of different color drawn in the first $n - (mt + m - 1)$ draws from an urn with balls of $m$ colors, initially containing $t + 1$ balls of each color, where we draw balls at random and replace each drawn ball together with a new ball of the same color (see, e.g., Johnson and Kotz [20], Section 4.5.1).

This distribution can be obtained by first taking a random vector $\mathbf{V}$ with the Dirichlet distribution above and then a multinomial variable with parameters $n - (mt + m - 1)$ and $\mathbf{V}$ ([20], Section 4.5.1). Using this representation, (3.6) follows immediately from the law of large numbers, even in the stronger form $\mathbf{V}_n/n \stackrel{\mathrm{a.s.}}{\longrightarrow} \mathbf{V}$.

It follows from (3.6) and Lemma 3.1 that

$$\frac{G_{V_{n,j}}(z)}{G_n(z)} \stackrel{\mathrm{d}}{\longrightarrow} V_j^{\lambda_1(z)-1}$$

if $z \in D_1$ and $E(z) \neq 0$. Hence, if $M_n(z)$ has a limit (in distribution) $Y(z)$ for some $z \in D_1$ with $E(z) \neq 0$, then this limit must satisfy the stochastic fixed point equation

$$Y(z) \stackrel{\mathrm{d}}{=} zV_1^{\lambda_1(z)-1}Y^{(1)}(z) + zV_2^{\lambda_1(z)-1}Y^{(2)}(z) + \cdots + zV_m^{\lambda_1(z)-1}Y^{(m)}(z),$$

(3.7)

where $Y^{(j)}(z)$ are independent copies of $Y(z)$ that are independent of $\mathbf{V}$. [Note that $z \in D_1$ and $E(z) \neq 0$ imply that $G_n(z) \to \infty$.]

In Section 9, we will show that this limit relation is actually true in a suitable domain, even in a strong sense, as asserted in Theorem 1.2. We will also see that we have a unique solution of this stochastic fixed point equation under the assumption $\mathbb{E}Y(z) = 1$ and a certain integrability condition.

## 4. Function spaces.

For functions defined on an interval $I \subseteq \mathbb{R}$, we use the space $D(I)$ of right-continuous functions with left-hand limits equipped with the Skorohod topology. A general definition of this topology is that $f_n \to f$ as $n \to \infty$ if and only if there exists a sequence $\lambda_n$ of strictly increasing continuous functions that map $I$ onto itself such that $\lambda_n(x) \to x$ and $f_n(\lambda_n(x)) \to f(x)$, uniformly on every compact subinterval of $I$; see, for example, [2], Chapter 3, $(I = [0,1])$, [24], [18], Chapter VI, [21], Appendix A2 $([0,\infty))$, [19], Section 2. It is of technical importance that this topology can be induced by a complete, separable metric [2], Chapter 14, [18], Theorem VI.1.14, [21], Theorem A2.2. Note that it matters significantly whether or not the endpoints are included in the interval $I$, but we can always reduce to the case of compact intervals because $f_n \to f$ in $D(I)$ if and only if $f_n \to f$ in $D(J_k)$ for an increasing sequence of compact intervals $J_k$ with $\bigcup J_k = I$. In particular, when $f$ is continuous, $f_n \to f$ in $D(I)$ if and only if $f_n \to f$ uniformly on every compact subinterval. Similarly, if $F_n$ and $F$ are random elements of $D(I)$ and $F$ is a.s. continuous, then $F_n \stackrel{\mathrm{d}}{\longrightarrow} F$ in $D(I)$ if and only if $F_n \stackrel{\mathrm{d}}{\longrightarrow} F$ in $D(J)$ for every compact subinterval $J \subseteq I$.

For analytic functions on a domain (i.e., a nonempty open connected set) $D \subseteq \mathbb{C}$, we will use two topological vector spaces.

- $\mathcal{H}(D)$ is the space of all analytic functions on $D$ with the topology of uniform convergence on compact sets. This topology can be defined by the family of seminorms $f \mapsto \sup_K |f|$, where $K$ ranges over the compact subsets of $D$. $\mathcal{H}(D)$ is a Fréchet space, that is, a locally convex space with a topology that can be defined by a complete metric, and it has (by Montel's theorem on normal families) the property that every closed bounded subset is compact (see, e.g., [28], Chapter 1.45, or [29], Example 10.II and Theorem 14.6). It is easily seen that the topology is separable [e.g., by regarding $\mathcal{H}(D)$ as a subspace of $C_0^\infty(D)$].

- $\mathcal{B}(D)$ is the *Bergman space* of all square-integrable analytic functions on $D$, equipped with the norm given by $\|f\|_{\mathcal{B}(D)}^2 = \int_D |f(z)|^2 \, dm(z)$, where $m$ is the two-dimensional Lebesgue measure. $\mathcal{B}(D)$ can be regarded as a closed subspace of $L^2(\mathbb{R}^2)$ and is thus a separable Hilbert space (see, e.g., [22], Chapter 1.4).

Since these spaces are metric spaces, we can use the general theory in, for example, Billingsley [2] or Kallenberg [21] for convergence in distribution of random functions in these spaces (equipped with their Borel $\sigma$-fields).

$\mathcal{B}(D)$ has the advantage of being a Hilbert space, which will be important for us later. On the other hand, $\mathcal{H}(D)$ is, in several ways, the natural space for analytic functions. One important technical advantage of $\mathcal{H}(D)$ is that it is easy to characterize tightness. Recall that a sequence $(W_n)$ of random variables in a metric space $\mathcal{S}$ is *tight* if for every $\varepsilon > 0$, there exists a compact subset $K \subseteq \mathcal{S}$ such that $\mathbb{P}(W_n \in K) > 1 - \varepsilon$ for every $n$. In a Polish space, that is, a complete separable metric space, tightness is equivalent to relative compactness (of the corresponding distributions) by Prohorov's theorem [2], Theorems 6.1 and 6.2, [21], Theorem 16.3. [Both $\mathcal{H}(D)$ and $\mathcal{B}(D)$ are Polish, by the properties above.]

**Lemma 4.1.** *Let $D$ be a domain in $\mathbb{C}$. A sequence $(W_n)$ of random analytic functions on $D$ is tight in $\mathcal{H}(D)$ if and only if the sequence $(\sup_{z \in K} |W_n(z)|)$ is tight for every compact $K \subset D$, that is, if and only if for every compact $K \subset D$ and every $\varepsilon > 0$, there exists an $M$ such that $\mathbb{P}(\sup_{z \in K} |W_n(z)| > M) < \varepsilon$ for all $n$.*

**Proof.** This is an easy consequence of the characterization of compact sets as closed bounded sets in $\mathcal{H}(D)$. We omit the details. $\square$

The embedding $\mathcal{B}(D) \to \mathcal{H}(D)$ is continuous [22], Lemma 1.4.1. Thus, convergence in distribution in $\mathcal{B}(D)$ implies convergence in $\mathcal{H}(D)$. Similarly, if $D' \subset D$ is a subdomain, then the restriction mappings $\mathcal{H}(D) \to \mathcal{H}(D')$ and $\mathcal{B}(D) \to \mathcal{B}(D')$ are continuous and thus convergence in distribution in $\mathcal{H}(D)$ or $\mathcal{B}(D)$ implies convergence (of the restrictions) in $\mathcal{H}(D')$ or $\mathcal{B}(D')$, respectively.

The following theorem is a converse, which makes it possible to reduce the proof of convergence in $\mathcal{H}$ to local arguments. In applications, it is convenient to let $D_x$ be a small disc with center $x$.

**Theorem 4.1.** *Let $D \subseteq \mathbb{C}$ be a domain. Suppose that $(W_n)$ is a sequence of random analytic functions on $D$ and that for each $x \in D$, there is an open subdomain $D_x$ with $x \in D_x \subset D$ and a random analytic function $Z_x$ on $D_x$ such that $W_n \xrightarrow{\mathrm{d}} Z_x$ in $\mathcal{H}(D_x)$ as $n \to \infty$. There then exists a random analytic function $Z$ on $D$ such that $W_n \xrightarrow{\mathrm{d}} Z$ in $\mathcal{H}(D)$ as $n \to \infty$ and the restriction $Z|_{D_x} \stackrel{\mathrm{d}}{=} Z_x$ for every $x$.*

To prove this, we use the following general measure-theoretic lemma, which we copy from Bousquet-Mélou and Janson [3], Lemma 7.1.

LEMMA 4.2. *Let $\mathcal{S}_1$ and $\mathcal{S}_2$ be two Polish spaces and let $\phi \colon \mathcal{S}_1 \to \mathcal{S}_2$ be an injective continuous map. If $(W_n)$ is a tight sequence of random elements of $\mathcal{S}_1$ such that $\phi(W_n) \xrightarrow{\mathrm{d}} Z$ in $\mathcal{S}_2$ for some random $Z \in \mathcal{S}_2$, then $W_n \xrightarrow{\mathrm{d}} W$ in $\mathcal{S}_1$ for some $W$ with $\phi(W) \stackrel{\mathrm{d}}{=} Z$.*

PROOF OF THEOREM 4.1. Let, for every $x \in D$, $D'_x$ be a small open disc with center $x$ such that $\overline{D'_x} \subset D_x$. Since $W_n \xrightarrow{\mathrm{d}} Z_x$ in $\mathcal{H}(D_x)$, the sequence $(W_n)$ is tight in $\mathcal{H}(D_x)$ for every $x \in D$. In particular, by Lemma 4.1, the sequence $(\sup_{D'_x} |W_n|)$ is tight.

If $K \subset D$ is compact, then $K$ can be covered by a finite number of the discs $D'_x$ and it follows that the sequence $(\sup_K |W_n|)$ is tight. Consequently, the sequence $(W_n)$ is tight in $\mathcal{H}(D)$ by Lemma 4.1.

We now fix $x \in D$ and apply Lemma 4.2 with $\mathcal{S}_1 = \mathcal{H}(D)$, $\mathcal{S}_2 = \mathcal{H}(D_x)$ and $\phi$ the restriction map. Note that $\phi$ is injective since the functions are analytic and $D$ is connected. The result follows. $\quad\square$

For future use, we include the following alternative characterization of completeness in an arbitrary complete metric space $\mathcal{S}$. If $A \subseteq \mathcal{S}$, let $A^\varepsilon$ denote the set $\{x \colon d(x, A) < \varepsilon\}$.

LEMMA 4.3. *If $\{W_\alpha\}$ is a family of random variables in a complete metric space $\mathcal{S}$, then the following are equivalent characterizations of tightness of $\{W_\alpha\}$:*

(i) *for every $\varepsilon > 0$, there exists a compact set $K \subseteq \mathcal{S}$ such that*

$$\mathbb{P}(X_\alpha \notin K) < \varepsilon \qquad \text{for every } \alpha;$$

(ii) *for every $\varepsilon, \delta > 0$, there exists a compact set $K \subseteq \mathcal{S}$ such that*

$$\mathbb{P}(X_\alpha \notin K^\delta) < \varepsilon \qquad \text{for every } \alpha;$$

(iii) *for every $\varepsilon, \delta > 0$, there exists a finite set $F \subseteq \mathcal{S}$ such that*

$$\mathbb{P}(X_\alpha \notin F^\delta) < \varepsilon \qquad \text{for every } \alpha.$$

PROOF. (i): This is the standard definition of tightness [2], Chapter 6.

(i) $\implies$ (ii): This is obvious.

(ii) $\implies$ (iii): Given $\varepsilon$ and $\delta$, let $K$ be as in (ii). Since $K$ is compact, there exists a finite set $F$ such that $K \subseteq F^\delta$ and thus $K^\delta \subseteq F^{2\delta}$. Hence, $\mathbb{P}(X_\alpha \notin F^{2\delta}) < \varepsilon$.

(iii) $\implies$ (i): Let $F_n$ be a finite set such that $\mathbb{P}(X_\alpha \notin F_n^{1/n}) < \varepsilon \cdot 2^{-n}$ for every $\alpha$, and let $K := \bigcap_{n \geq 1} F_n^{1/n}$. $K$ is then closed and totally bounded, thus compact, and $\mathbb{P}(X_\alpha \notin K) \leq \sum_n \mathbb{P}(X_\alpha \notin F_n^{1/n}) < \varepsilon$ for every $\alpha$. $\quad\square$

**5. The Zolotarev metric on a Hilbert space.** We recall the definition of the Zolotarev metric for probability measures in a Banach space; see Zolotarev [30].

If $B$ and $B_1$ are Banach spaces and $f : U \to B_1$ is a function defined on an open subset $U \subseteq B$, then $f$ is said to be (Fréchet) differentiable at a point $x \in U$ if there exists a linear operator $Df(x) : B \to B_1$ such that $\|f(x+y) - f(x) - Df(x)y\|_{B_1} = o(\|y\|_B)$ as $\|y\|_B \to 0$. Further, $f$ is differentiable in $U$ if it is differentiable for every $x \in U$. $Df$ is then a function $U \to L(B, B_1)$ and we may talk about its derivative $D^2f = DDf$, and so on. Note that the $m$th derivative $D^m f$ (if it exists) is a function from $U$ into the space of multilinear mappings $B^m \to B_1$. Let $C^m(B, B_1)$ denote the space of $m$ times continuously differentiable functions $f : B \to B_1$.

Given a Banach space $B$ and a real number $s > 0$, write $s = m + \alpha$ with $0 < \alpha \le 1$ and $m := \lceil s \rceil - 1 \in \mathbb{N}_{\ge 0}$, and define

$$\mathcal{F}_s := \{ f \in C^m(B, \mathbb{R}) : \|D^m f(x) - D^m f(y)\| \le \|x - y\|^{\alpha}, x, y \in B \}.$$

We will also write $\mathcal{F}_s^* := \{ f \in C^m(B, \mathbb{R}) : cf \in \mathcal{F}_s \text{ for some } c > 0 \}$.

The Zolotarev metric $\zeta_s$ is a distance between distributions, but it is often convenient to talk about it as a distance between random variables, keeping in mind that only their distributions matter. For two random variables $X$ and $Y$ with values in $B$, or for their corresponding distributions $\mathcal{L}(X)$ and $\mathcal{L}(Y)$, the Zolotarev metric $\zeta_s$ is defined by

$$(5.1) \qquad \zeta_s(X, Y) := \zeta_s(\mathcal{L}(X), \mathcal{L}(Y)) := \sup_{f \in \mathcal{F}_s} |\mathbb{E}(f(X) - f(Y))|.$$

Note that this distance may be infinite, but it is easily seen, by a Taylor expansion of $f$, that it is finite if $\mathbb{E}\|X\|^s < \infty$, $\mathbb{E}\|Y\|^s < \infty$, and $X$ and $Y$ have the same moments up to order $m$, where the $k$th moment of $X$ is $\mathbb{E}X^{\otimes k}$, regarded as an element of the $k$th (completed) projective tensor power $B^{\otimes k}$.

REMARK 5.1. The dual space of $B^{\otimes k}$ is the space of bounded multilinear mappings $B^k \to \mathbb{R}$. Hence, $\mathbb{E}X^{\otimes k} = \mathbb{E}Y^{\otimes k}$ if and only if $\mathbb{E}g(X, \ldots, X) = \mathbb{E}g(Y, \ldots, Y)$ for every bounded multilinear mapping $B^k \to \mathbb{R}$. See, for example, [29], Chapter 45 for facts on tensor products.

We define, for a given sequence $\mathbf{z} = (z_1, \ldots, z_m)$ with $z_k \in B^{\otimes k}$, $k = 1, \ldots, m$,

$$\mathcal{P}_{s,\mathbf{z}}(B) := \{ \mathcal{L}(X) : \mathbb{E}\|X\|^s < \infty, \mathbb{E}X^{\otimes k} = z_k, k = 1, \ldots, m \},$$

that is, the set of probability measures on $B$ with finite absolute $s$th moment and moments $z_1, \ldots, z_k$. Thus, $\zeta_s$ is finite on each $\mathcal{P}_{s,\mathbf{z}}(B)$ and it is obviously a semi-metric there.

We are mainly concerned with the cases $0 < s \leq 1$ when $m = 0$ and $1 < s \leq 2$ when $m = 1$. In these cases, we write $\mathcal{P}_s(B)$ and $\mathcal{P}_{s,z}(B)$, respectively, where $z \in H$ is the mean.

For a general Banach space $B$, we do not know whether $\zeta_s$ is always a complete metric on $\mathcal{P}_{s,\mathbf{z}}(B)$. Moreover, according to Bentkus and Rachkauskas [1], it is not hard to show that in a general Banach space, convergence in $\zeta_s$ does not imply weak convergence (convergence in distribution) when $s > 1$, although we do not know of any explicit counterexample. (It is easy to see that convergence in $\zeta_s$ for $0 < s \leq 1$ implies weak convergence, by the proof of Theorem 2.1(ii) $\Longrightarrow$ (iii) in [2].) We will therefore, in the sequel, restrict ourselves to separable Hilbert spaces, where we can show these desirable properties.

THEOREM 5.1. *If $H$ is a separable Hilbert space and $s > 0$, then $\zeta_s$ is a complete metric on the set $\mathcal{P}_{s,\mathbf{z}}(H)$ of all probability measures on $H$ with a finite $s$th absolute moment and given $k$th moments $z_k$, $1 \leq k < s$. Moreover, if $X_n, X$ are $H$-valued random variables with distributions in $\mathcal{P}_{s,\mathbf{z}}(H)$ and $\zeta_s(X_n, X) \to 0$, then $X_n \overset{\mathrm{d}}{\longrightarrow} X$.*

The final assertion is proved by Giné and León [16]. For completeness, we include a short proof using lemmas needed for the first part.

PROOF OF THEOREM 5.1. First, note that $\zeta_s$ is a metric on $\mathcal{P}_{s,\mathbf{z}}(H)$ [30]; the fact that $\zeta_s(\mu, \nu) = 0$ implies $\mu = \nu$ for probability measures $\mu, \nu \in \mathcal{P}_{s,\mathbf{z}}(H)$ is well known and follows easily because $x \mapsto e^{\mathrm{i}\langle x,y\rangle} \in \mathcal{F}_s^*$ for every $y \in H$, thus if $\zeta_s(\mu, \nu) = 0$, then, by (5.1), the characteristic functions $\int e^{\mathrm{i}\langle x,y\rangle} d\mu(x)$ and $\int e^{\mathrm{i}\langle x,y\rangle} d\nu(x)$ are equal, which implies that all finite-dimensional projections coincide for $\mu$ and $\nu$, and $\mu = \nu$ then follows by a monotone class argument (see, e.g., [23], Section 2.1).

We continue by constructing some other functions in $\mathcal{F}_s^*$. Taking small positive multiples of them, we thus obtain functions in $\mathcal{F}_s$.

LEMMA 5.1. *Let $\varphi \in C^\infty(\mathbb{R})$ with $\varphi(x) = x^2$ for $|x| < 1/2$, $\varphi(x) = 1$ for $|x| > 1$ and $1/4 \leq \varphi(x) \leq 1$ for $1/2 \leq |x| \leq 1$. Then $f(x) := \varphi(\|x\|) \in \mathcal{F}_s^*$ for every $s > 0$.*

PROOF. First, note that $x \mapsto \|x\|^2$ is infinitely differentiable on $H$. (In fact, the third derivative vanishes.) Hence, if $g \colon \mathbb{R} \to \mathbb{R}$ is any $C^\infty$ function, then $g(\|x\|) = g((\|x\|^2)^{1/2})$ is infinitely differentiable on $H \setminus \{0\}$.

Consequently, $\varphi(\|x\|)$ is infinitely differentiable both in $\{x \colon \|x\| < 1/2\}$ and in $\{x \colon \|x\| > 0\}$, and thus everywhere. Further, any derivative of order $\geq 1$ vanishes for $\|x\| > 1$ and is bounded on $\|x\| \leq 1$, hence it is globally

bounded. In particular, $D^m f$ is both bounded and has a bounded derivative, which implies that $f \in \mathcal{F}_s^*$. (Consider the cases $\|x - y\| \leq 1$ and $\|x - y\| > 1$ separately.) $\square$

LEMMA 5.2. *Let $\psi \in C^\infty(\mathbb{R})$ with $\psi(x) \geq 0$, $\psi(x) = 0$ for $|x| \leq 1/2$ and $\varphi(x) = |x|^s$ for $|x| \geq 1$. Then $f(x) := \psi(\|x\|) \in \mathcal{F}_s^*$ for every $s > 0$.*

PROOF. It is easily seen, as in Lemma 5.1, that $f$ is infinitely differentiable and that $D^k f(x) = O(\|x\|^{s-k})$ for every fixed $k \geq 0$. Hence, if $x, y \in H$ with $\|x - y\| < \|x\|/2$, then

$$D^{m+1} f(z) = O(\|z\|^{s-m-1}) = O(\|x\|^{\alpha - 1})$$

for $z \in [x, y]$ and thus

$$\|D^m f(x) - D^m f(y)\| = O(\|x - y\|\|x\|^{\alpha - 1}) = O(\|x - y\|^\alpha).$$

The same holds by symmetry if $\|x - y\| < \|y\|/2$.

Finally, if $\|x - y\| \geq \frac{1}{2}\|x\|, \frac{1}{2}\|y\|$, then

$$\begin{aligned}
\|D^m f(x) - D^m f(y)\| &\leq \|D^m f(x)\| + \|D^m f(y)\| \\
&= O(\|x\|^{s-m} + \|y\|^{s-m}) \\
&= O(\|x - y\|^\alpha).
\end{aligned}$$

Thus, $f \in \mathcal{F}_s^*$. $\square$

In the following Lemmas 5.3–5.6, we assume that $\{\mu_n\}_1^\infty$ is a sequence of probability measures in $\mathcal{P}_{s,\mathbf{z}}(H)$ that is a Cauchy sequence for $\zeta_s$ and let $\{X_n\}_1^\infty$ be $H$-valued random variables such that $X_n$ has the distribution $\mu_n$.

LEMMA 5.3. *The random variables $\|X_n\|^s$ are uniformly integrable.*

PROOF. Let $f(x) = \psi(\|x\|)$ be as in Lemma 5.2 (for some fixed choice of $\psi$); by Lemma 5.2, $cf \in \mathcal{F}_s$ for some constant $c > 0$. Let, for $r > 0$, $f_r(x) = cr^s f(x/r)$. Then, as is easily seen, $f_r \in \mathcal{F}_s$ for every $r > 0$. Further, $f_r(x) = c\|x\|^s$ when $\|x\| > r$, but $f_r(x) = 0$ when $\|x\| < r/2$.

Since $|f_r(X_m)| \leq cC\|X_m\|^s$ for all $r$, where $C := \sup_{x \neq 0} f(x)/\|x\|^s < \infty$, and $f_r(X_m) \to 0$ as $r \to \infty$, dominated convergence yields $\mathbb{E} f_r(X_m) \to 0$ as $r \to \infty$ for every fixed $m$.

Now, let $\varepsilon > 0$ and choose $N$ such that $\zeta_s(X_n, X_N) < \varepsilon$ for $n \geq N$. Choose $r$ such that $\mathbb{E} f_r(X_m) < \varepsilon$ for $m = 1, \ldots, N$. Then, for every $n > N$,

$$\mathbb{E} f_r(X_n) \leq \mathbb{E} f_r(X_N) + \zeta_s(X_n, X_N) < 2\varepsilon.$$

Thus, $\mathbb{E} f_r(X_n) < 2\varepsilon$ for all $n$. Consequently, for all $n$,

$$\mathbb{E}(\|X_n\|^s \mathbf{1}[\|X_n\| > r]) \leq c^{-1} \mathbb{E} f_r(X_n) < 2c^{-1}\varepsilon.$$

Since $\varepsilon > 0$ is arbitrary, this shows the desired uniform integrability. $\square$

Lemma 5.4. *The sequence $(X_n)$ is tight.*

Proof. Let $f(x)$ be as in Lemma 5.1, and let $c > 0$ be such that $cf \in \mathcal{F}_s$. Let $\varepsilon, \eta > 0$ and assume that $\varepsilon, \eta < 1/2$. Let $\delta := \varepsilon^2 \eta / 6 < 1/2$.

Choose $N$ such that $\zeta_s(X_n, X_N) < c\delta$ for $n \geq N$. Since $H$ is complete and separable, each $X_n$ is tight [2], Theorem 1.4. Hence, there exists, by Lemma 4.3, for each $n$, a finite set $F_n \subset H$ such that $\mathbb{P}(X_n \notin F_n^\delta) < \delta$. Let $F := \bigcup_1^N F_n$.

Let $M$ be the subspace of $H$ spanned by $F$, let $P$ be the orthogonal projection onto $M$ and let $Q = I - P$ be the complementary projection. If $X_n \in F^\delta$, then $\|Q(X_n)\| < \delta$ and thus $f(Q(X_n)) = \|Q(X_n)\|^2 < \delta^2$. Hence,

$$(5.2) \qquad \mathbb{E} f(Q(X_n)) \leq \delta^2 + \mathbb{P}(X_n \notin F^\delta) < 2\delta, \qquad n \leq N.$$

Since $cf \in \mathcal{F}_s$ and $Q$ is a linear operator with norm at most 1, it is easily seen that $cf \circ Q \in \mathcal{F}_s$ too. Hence, for $n \geq N$,

$$|\mathbb{E} f(Q(X_n)) - \mathbb{E} f(Q(X_N))| \leq c^{-1} \zeta_s(X_n, X_N) < \delta$$

and, by (5.2),

$$\mathbb{E} f(Q(X_n)) \leq \mathbb{E} f(Q(X_N)) + \delta \leq 3\delta, \qquad n \geq N.$$

Combining this with (5.2), we see that $\mathbb{E} f(Q(X_n)) \leq 3\delta$ for all $n$. Hence, by Markov's inequality,

$$(5.3) \qquad \mathbb{P}(\|Q(X_n)\| > \varepsilon) \leq \mathbb{P}(f(Q(X_n)) > \varepsilon^2) \leq 3\delta/\varepsilon^2 = \eta/2.$$

Next, the random variables $P(X_n)$ lie in the finite-dimensional space $M$ and $\sup_n \mathbb{E}\|P(X_n)\|^s \leq \sup_n \mathbb{E}\|X_n\|^s < \infty$, by Lemma 5.3. Let $K_R := \{x \in M : \|x\| \leq R\}$. Then $K_R$ is compact and it follows from Markov's inequality that if $R$ is large enough, then $\mathbb{P}(P(X_n) \notin K_R) < \eta/2$ for every $n$ and consequently, recalling (5.3),

$$(5.4) \qquad \mathbb{P}(X_n \notin K_R^\varepsilon) \leq \mathbb{P}(P(X_n) \notin K_R) + \mathbb{P}(\|Q(X_n)\| \geq \varepsilon) < \eta.$$

We have shown that for every $\varepsilon, \eta > 0$, there exists a compact set $K_R$ such that (5.4) holds for all $n$. (Clearly, we may assume that $\varepsilon, \eta < 1/2$, as we have done.) By Lemma 4.3, the family $\{X_n\}$ is tight. □

Lemma 5.5. *If, further, $X_n \overset{d}{\longrightarrow} X$ for some $H$-valued random variable $X$, then $\mathcal{L}(X) \in \mathcal{P}_{s,\mathbf{z}}(H)$ and $\mathbb{E} f(X_n) \to \mathbb{E} f(X)$ for every $f \in \mathcal{F}_s$.*

Proof. First, by Fatou's lemma and Lemma 5.3,

$$\mathbb{E}\|X\|^s \leq \liminf_{n \to \infty} \mathbb{E}\|X_n\|^s \leq \sup_n \mathbb{E}\|X_n\|^s < \infty.$$

Next, if $f \in \mathcal{F}_s$, then $f(X_n) \xrightarrow{\mathrm{d}} f(X)$ because $f$ is continuous. Further, a Taylor expansion yields (see [30], (18)) $|f(x)| \leq K + K\|x\|^s$ for some $K$ (depending on $f$). Hence, $|f(X_n)| \leq K + K\|X_n\|^s$ and $\{f(X_n)\}$ is uniformly integrable by Lemma 5.3. Consequently, $\mathbb{E}f(X_n) \to \mathbb{E}f(X)$.

If $f(x) = g(x, \dots, x)$ for some continuous multilinear mapping $g : B^k \to \mathbb{R}$ with $k \leq m$, then, for every $n$, using the duality between tensor powers and multilinear mappings in Remark 5.1,

$$\mathbb{E}f(X_n) = \mathbb{E}\langle g, X_n^{\otimes k}\rangle = \langle g, \mathbb{E}X_n^{\otimes k}\rangle = \langle g, z_k\rangle.$$

Further, $D^m f$ is constant and thus $f \in \mathcal{F}_s$, consequently,

$$\langle g, \mathbb{E}X^{\otimes k}\rangle = \mathbb{E}f(X) = \lim_n \mathbb{E}f(X_n) = \langle g, z_k\rangle.$$

Since $g$ is arbitrary, this implies that $\mathbb{E}X^{\otimes k} = z_k$ and thus $\mathcal{L}(X) \in \mathcal{P}_{s,\mathbf{z}}(H)$. $\square$

LEMMA 5.6. *If $X_n \xrightarrow{\mathrm{d}} X$ for some $H$-valued random variable $X$, then $\zeta_s(X_n, X) \to 0$.*

PROOF. Let $\varepsilon > 0$ and choose $N$ such that $\zeta_s(X_n, X_m) < \varepsilon$ if $n, m \geq N$. For any $f \in \mathcal{F}_s$ and any $n, m \geq N$, we thus have $|\mathbb{E}f(X_n) - \mathbb{E}f(X_m)| < \varepsilon$. Letting $m \to \infty$, we thus obtain, by Lemma 5.5, $|\mathbb{E}f(X_n) - \mathbb{E}f(X)| \leq \varepsilon$ for $n \geq N$ and every $f \in \mathcal{F}_s$. Thus, $\zeta_s(X_n, X) \leq \varepsilon$ for $n \geq N$. $\square$

We may now complete the proof of Theorem 5.1. First, assume that $(\mu_n)$ is a Cauchy sequence in $\mathcal{P}_{s,\mathbf{z}}(H)$, equipped with metric $\zeta_s$. Let $X_n$ be random variables with the distributions $\mu_n$. By Lemma 5.4, the sequence $(X_n)$ is tight, so, by Prohorov's theorem, there exists a subsequence that converges in distribution to some $H$-valued random variable $X$. Let $\mu$ be the distribution of $X$. Considering this subsequence only, we see by Lemma 5.5 that $\mu \in \mathcal{P}_{s,\mathbf{z}}(H)$ and by Lemma 5.6 that $\mu_n \to \mu$ in $\mathcal{P}_{s,\mathbf{z}}(H)$ along the subsequence. Since $(\mu_n)$ is a Cauchy sequence, the full sequence also converges. Hence $\mathcal{P}_{s,\mathbf{z}}(H)$ is complete.

Second, assume that $X_n$ and $X$ are $H$-valued random variables with distributions in $\mathcal{P}_{s,\mathbf{z}}(H)$ such that $\zeta_s(X_n, X) \to 0$. In particular, the distributions $\mathcal{L}(X_n)$ form a Cauchy sequence in $\mathcal{P}_{s,\mathbf{z}}(H)$, so, by Lemma 5.4, the sequence is tight. If a subsequence converges in distribution to some random variable $Y$, then Lemma 5.6 shows that $\zeta_s(X_n, Y) \to 0$ along the subsequence and thus $\zeta_s(X, Y) = 0$, so $Y \xrightarrow{\mathrm{d}} X$. Hence, all subsequence limits of $(X_n)$ have the distribution of $X$ and since the sequence is tight, this means that $X_n \xrightarrow{\mathrm{d}} X$. $\square$

We will later use an upper bound of $\zeta_s$ by the minimal $L^s$-metric $\ell_s$; see Zolotarev [30] for similar bounds. The $\ell_s$ metric, $s > 0$, is defined on all random variables $X$ and $Y$ with values in $B$ and satisfying $\mathbb{E}\|X\|^s$, $\mathbb{E}\|Y\|^s < \infty$ by

$$\ell_s(X, Y) := \ell_s(\mathcal{L}(X), \mathcal{L}(Y))$$
$$:= \inf\{(\mathbb{E}\|X' - Y'\|^s)^{(1/s)\wedge 1} : \mathcal{L}(X') = \mathcal{L}(X), \mathcal{L}(Y') = \mathcal{L}(Y)\}.$$

LEMMA 5.7.   *For all* $\mathcal{L}(X), \mathcal{L}(Y) \in \mathcal{P}_{s,\mathbf{z}}(B)$ *and* $s > 1$, *we have*

$$(5.5) \qquad \zeta_s(X, Y) \leq ((\mathbb{E}\|X\|^s)^{1-1/s} + (\mathbb{E}\|Y\|^s)^{1-1/s})\ell_s(X, Y).$$

*For* $0 < s \leq 1$, *we have*

$$\zeta_s(X, Y) \leq \ell_s(X, Y).$$

PROOF.   For $s > 1$ and arbitrary $f \in \mathcal{F}_s$, we define

$$g(x) := f(x) - f(0) - Df(0)(x) - \cdots - \frac{1}{m!}D^m f(0)(x, \ldots, x).$$

Thus, we have

$$g(0) = Dg(0) = \cdots = D^m g(0) = 0, \qquad \|D^m g(x) - D^m g(y)\| \leq \|x - y\|^\alpha.$$

This implies, by backward induction on $k$, that

$$\|D^k g(x)\| \leq \|x\|^{s-k}, \qquad 0 \leq k \leq m.$$

Thus, with $Z := Y - X$, we obtain, for an appropriate $0 \leq \vartheta \leq 1$,

$$|g(Y) - g(X)| = |g(X + Z) - g(X)| = |Dg(X + \vartheta Z)(Z)|$$
$$\leq \|Dg(X + \vartheta Z)\|\|Z\| \leq \|X + \vartheta Z\|^{s-1}\|Z\|$$
$$\leq (\|X\|^{s-1} + \|Y\|^{s-1})\|Z\|.$$

This implies, using Hölder's inequality, that

$$|\mathbb{E}[f(Y) - f(X)]| = |\mathbb{E}[g(Y) - g(X)]| \leq \mathbb{E}[(\|X\|^{s-1} + \|Y\|^{s-1})\|Z\|]$$
$$\leq ((\mathbb{E}\|X\|^s)^{1-1/s} + (\mathbb{E}\|Y\|^s)^{1-1/s})(\mathbb{E}\|Z\|^s)^{1/s}.$$

Taking the supremum over all $f \in \mathcal{F}_s$ and the infimum over all realizations of $\mathcal{L}(X)$ and $\mathcal{L}(Y)$, we obtain (5.5).

For $0 < s \leq 1$, we have $|\mathbb{E}[f(Y) - f(X)]| \leq \mathbb{E}|f(Y) - f(X)| \leq \|X - Y\|^s$. This implies that $\zeta_s(X, Y) \leq \ell_s(X, Y)$.   $\square$

**6. Contraction method for Hilbert spaces.** In this section, we extend the contraction method as developed for the Zolotarev metric on $\mathbb{R}^d$ in [27] to random variables in a separable Hilbert space $H$. We denote by $\mathcal{P}(H)$ the set of all probability distributions on $H$. The limit distributions occurring subsequently are characterized as fixed points of the maps

(6.1)
$$T : \mathcal{P}(H) \to \mathcal{P}(H),$$

$$\mathcal{L}(Z) \mapsto \mathcal{L}\left( \sum_{r=1}^m A_r^*(Z^{(r)}) + b^* \right),$$

where $m \geq 1$ is an integer, $A_1^*, \ldots, A_m^*$ are random linear operators in $H$, $b^*$ is a random variable in $H$ and $(A_1^*, \ldots, A_m^*, b^*)$, $Z^{(1)}, \ldots, Z^{(m)}$ are independent, with $\mathcal{L}(Z^{(r)}) = \mathcal{L}(Z)$ for $r = 1, \ldots, m$.

We write $\|A\|_{\mathrm{op}} := \sup_{\|x\|=1} \|Ax\|$ for a linear operator $A$ in $H$. We say that $A_r^*$ is $s$-*integrable* if $\mathbb{E}\|A_r^*\|_{\mathrm{op}}^s < \infty$. Furthermore, we abbreviate $\mathcal{P}_s := \mathcal{P}_s(H)$ and $\mathcal{P}_{s,0} := \mathcal{P}_{s,0}(H)$, as defined in Section 5.

LEMMA 6.1. *Let* $(A_1^*, \ldots, A_m^*, b^*)$ *be as in* (6.1) *and $s$-integrable for some* $0 < s \leq 2$. *For* $0 < s \leq 1$, *we have* $T(\mathcal{P}_s) \subseteq \mathcal{P}_s$. *For* $1 < s \leq 2$ *and* $\mathbb{E}b^* = 0$, *we have* $T(\mathcal{P}_{s,0}) \subseteq \mathcal{P}_{s,0}$.

PROOF. The existence of moments of order $s$ of $T(\mathcal{L}(Z))$ with $\mathbb{E}\|Z\|^s$ follows by the independence of $A_r^*$ and $Z^{(r)}$. For $1 < s \leq 2$ and $\mathbb{E}b^* = 0$, we obtain that $T(\mathcal{L}(Z))$ is centered. $\square$

LEMMA 6.2. *Let* $(A_1^*, \ldots, A_m^*, b^*)$ *be as in* (6.1) *and $s$-integrable for some* $0 < s \leq 2$. *Assume that*

(6.2)
$$\mathbb{E} \sum_{r=1}^m \|A_r^*\|_{\mathrm{op}}^s < 1.$$

*If* $0 < s \leq 1$, *then the restriction of $T$ to $\mathcal{P}_s$ is a strict contraction. If* $1 < s \leq 2$ *and* $\mathbb{E}b^* = 0$, *then the restriction of $T$ to $\mathcal{P}_{s,0}$ is a strict contraction.*

PROOF. This is similar to the proof of Lemma 3.1 in [27]. Note that for a linear operator $A$ in $H$ and $\mathcal{L}(X), \mathcal{L}(Y) \in \mathcal{P}_{s,\mathbf{z}}$, we have $\zeta_s(A(X), A(Y)) \leq \|A\|_{\mathrm{op}}^s \zeta_s(X, Y)$ (cf. Zolotarev [31], Theorem 3). $\square$

Lemma 6.2 and Theorem 5.1 imply that the restrictions of $T$ in Lemma 6.2 have unique fixed points in $\mathcal{P}_s$ and $\mathcal{P}_{s,0}$, respectively.

We consider a sequence $(X_n)_{n \geq 0}$ of random variables in $H$ satisfying the recurrence

(6.3)
$$X_n \stackrel{d}{=} \sum_{r=1}^m A_r^{(n)}\big(X_{I_r^{(n)}}^{(r)}\big) + b^{(n)}, \qquad n \geq n_0,$$

where $n_0 \geq 1$, $A_r^{(n)}$ are random linear operators in $H$, $b^{(n)}$ is a random variable in $H$, $I^{(n)} = (I_1^{(n)}, \ldots, I_m^{(n)})$ is a vector of random integers with $I_r^{(n)} \in \{0, \ldots, n\}$, we have that $(X_j^{(1)}), \ldots, (X_j^{(m)})$, $(A_1^{(n)}, \ldots, A_m^{(n)}, b^{(n)}, I^{(n)})$ are independent and $\mathcal{L}(X_j^{(r)}) = \mathcal{L}(X_j)$ for all $r$ and $j$. We then have the following extension of Theorem 4.1 in [27].

THEOREM 6.1. *Let $(X_n)$ be as in* (6.3) *with all quantities there being $s$-integrable for some $0 < s \leq 2$. For $1 < s \leq 2$, assume that $\mathbb{E}X_n = 0$ for all $n \geq 0$. Assume that, for appropriately chosen $(A_1^*, \ldots, A_m^*, b^*)$, we have*

$$(6.4) \qquad \mathbb{E}\|A_r^{(n)} - A_r^*\|_{op}^s \to 0, \qquad \mathbb{E}\|b^{(n)} - b^*\|^s \to 0,$$

$$(6.5) \qquad \mathbb{E}\sum_{r=1}^m \|A_r^*\|_{op}^s < 1,$$

$$(6.6) \qquad \mathbb{E}[\mathbf{1}_{\{I_r^{(n)} \leq \ell\}}\|A_r^{(n)}\|_{op}^s] \to 0,$$

*for all $\ell \in \mathbb{N}$ and $r = 1, \ldots, m$. We then have*

$$(6.7) \qquad \zeta_s(X_n, X) \to 0, \qquad n \to \infty,$$

*where $\mathcal{L}(X)$ is the unique fixed point of $T$ in $\mathcal{P}_s$ for $0 < s \leq 1$ and in $\mathcal{P}_{s,0}$ for $1 < s \leq 2$.*

PROOF. The proof of Theorem 4.1 in [27] can be directly extended. For $1 < s \leq 2$, from $\mathbb{E}X_n = 0$ and (6.3), we obtain $\mathbb{E}b^{(n)} = 0$, thus (6.4) implies $\mathbb{E}b^* = 0$. Therefore, Lemma 6.2 implies existence and uniqueness of the fixed point of $T$, as claimed in the theorem for all $0 < s \leq 2$.

We introduce the accompanying sequence

$$(6.8) \qquad Q_n := \stackrel{d}{=} \sum_{r=1}^m A_r^{(n)}(X^{(r)}) + b^{(n)}, \qquad n \geq n_0,$$

where $(A_1^{(n)}, \ldots, A_m^{(n)}, b^{(n)}, I^{(n)})$, $X^{(1)}, \ldots, X^{(m)}$, $(X_n^{(1)}), \ldots, (X_n^{(m)})$ are independent, with $\mathcal{L}(X^{(r)}) = \mathcal{L}(X)$ for $r = 1, \ldots, m$. We obtain that $\mathbb{E}\|Q_n\|^s < \infty$ and, for $1 < s \leq 2$, that $\mathbb{E}Q_n = 0$. Thus, the $\zeta_s$ distances between $X_n, Q_n$ and $X$ are finite for $n \geq n_0$. We obtain, from the triangle inequality, that

$$(6.9) \qquad \zeta_s(X_n, X) \leq \zeta_s(X_n, Q_n) + \zeta_s(Q_n, X).$$

First, we show that $\zeta_s(Q_n, X) \to 0$. For this, note that we have $\mathbb{E}\|X\|^s < \infty$ and $\sup_{n \geq n_0} \mathbb{E}\|Q_n\|^s < \infty$ by representation (6.8), independence and (6.4). Hence, Lemma 5.7 implies that $\zeta_s(Q_n, X) \leq C\ell_s(Q_n, X)$, with some constant $C$ that is independent of $n$. For a random variable $Y$ in $H$, we denote the $L^s$-norm of $Y$ by $\|Y\|_s := (\mathbb{E}\|Y\|^s)^{1/s}$. Using (6.4), we obtain from $X \overset{d}{=} \sum A_r^* X^{(r)} + b^*$ that

$$
\begin{aligned}
\ell_s(Q_n, X) &\leq \left\| \sum_{r=1}^{m} (A_r^* - A_r^{(n)})(X^{(r)}) \right\|_s^{s \wedge 1} + \|b^{(n)} - b^*\|_s^{s \wedge 1} \\
&\leq \sum_{r=1}^{m} \|\|A_r^* - A_r^{(n)}\|_{\mathrm{op}} \|X\|\|_s^{s \wedge 1} + \|b^{(n)} - b^*\|_s^{s \wedge 1} \\
&\to 0, \qquad n \to \infty.
\end{aligned}
$$

Next, we bound the first summand in (6.9). Let $\Upsilon_n$ denote the joint distribution of $(A_1^{(n)}, \ldots, A_m^{(n)}, b^{(n)}, I^{(n)})$ and write $\alpha = (\alpha_1, \ldots, \alpha_m)$, $j = (j_1, \ldots, j_m)$. We then obtain, for $n \geq n_0$,

$$
\begin{aligned}
&\zeta_s(X_n, Q_n) \\
&= \zeta_s\left( \sum_{r=1}^{m} A_r^{(n)}(X_{I_r^{(n)}}^{(r)}) + b^{(n)}, \sum_{r=1}^{m} A_r^{(n)}(X^{(r)}) + b^{(n)} \right) \\
(6.10) \qquad &\leq \int \zeta_s\left( \sum_{r=1}^{m} \alpha_r(X_{j_r}^{(r)}), \sum_{r=1}^{m} \alpha_r(X^{(r)}) \right) d\Upsilon_n(\alpha, \beta, j) \\
&\leq \int \sum_{r=1}^{m} \|\alpha_r\|_{\mathrm{op}}^s \zeta_s(X_{j_r}, X) \, d\Upsilon_n(\alpha, \beta, j) \\
&\leq \left( \mathbb{E} \sum_{r=1}^{m} \|A_r^{(n)}\|_{\mathrm{op}}^s \right) \max_{0 \leq j \leq n} \zeta_s(X_j, X).
\end{aligned}
$$

Thus, using (6.9), it follows that

$$
\zeta_s(X_n, X) \leq \left( \mathbb{E} \sum_{r=1}^{m} \|A_r^{(n)}\|_{\mathrm{op}}^s \right) \max_{0 \leq j \leq n} \zeta_s(X_j, X) + o(1).
$$

We let $d_n := \zeta_s(X_n, X)$. We show that the sequence $(d_n)$ is bounded. With $r_n := \zeta_s(Q_n, X)$, $\xi_n := \mathbb{E}\sum_{r=1}^{m} \|A_r^{(n)}\|_{\mathrm{op}}^s$, $d_n^* := \max_{0 \leq j \leq n} d_j$ and $R := \sup_{n \geq n_0} r_n$, we obtain from (6.9) and (6.10) that

$$
(6.11) \qquad d_n \leq \xi_n d_n^* + r_n, \qquad n \geq n_0.
$$

By (6.4) and (6.5), there exist $\xi < 1$ and $n_1 \geq n_0$ such that $\xi_n \leq \xi$ for all $n \geq n_1$. Let $B := d_{n_1}^* + R/(1 - \xi)$. We claim that $d_n^* \leq B$ for all $n \geq 0$. If this fails to hold, then for some $n$, we have $d_{n-1}^* \leq B < d_n^* = d_n \vee d_{n-1}^*$. Hence, $d_n = d_n^* > B$. Moreover, we have $n > n_1$. Hence, (6.11) yields $d_n \leq \xi d_n + R$, thus $d_n \leq R/(1 - \xi) \leq B$, contradicting our assumption. Consequently, we have $d_n^* \leq B$ for all $n \geq 0$.

Let $\eta := \limsup_{n\to\infty} d_n$ and let $\varepsilon > 0$ be arbitrary. There exists an $\ell \in \mathbb{N}$ with $d_n \le \eta + \varepsilon$ for all $n \ge \ell$. We deduce, using (6.10), (6.9), (6.6) and (6.4), that

$$
\begin{aligned}
d_n &\le \int \sum_{r=1}^{m} \mathbf{1}_{\{0 \le j_r \le \ell\}} \|\alpha_r\|_{\mathrm{op}}^s \zeta_s(X_{j_r}, X) \, d\Upsilon_n(\alpha, \beta, j) \\
&\quad + \int \sum_{r=1}^{m} \mathbf{1}_{\{\ell < j_r \le n\}} \|\alpha_r\|_{\mathrm{op}}^s \zeta_s(X_{j_r}, X) \, d\Upsilon_n(\alpha, \beta, j) + r_n \\
&\le \mathbb{E} \sum_{r=1}^{m} \left(\mathbf{1}_{\{I_r^{(n)} \le \ell\}} \|A_r^{(n)}\|_{\mathrm{op}}^s\right) + (\eta + \varepsilon) \mathbb{E} \sum_{r=1}^{m} \|A_r^{(n)}\|_{\mathrm{op}}^s + r_n \\
&\le \left(\mathbb{E} \sum_{r=1}^{m} \|A_r^*\|_{\mathrm{op}}^s\right)(\eta + \varepsilon) + o(1).
\end{aligned}
$$

With $n \to \infty$, we obtain

$$
\eta \le \left(\mathbb{E} \sum_{r=1}^{m} \|A_r^*\|_{\mathrm{op}}^s\right)(\eta + \varepsilon).
$$

Since $\varepsilon > 0$ is arbitrary and $\mathbb{E} \sum_{r=1}^{m} \|A_r^*\|_{\mathrm{op}}^s < 1$, we obtain $\eta = 0$. Hence, $\zeta_s(X_n, X) \to 0$. $\square$

REMARK 6.1.  Note that the conditions (6.4) on the $L^s$-convergence of the coefficients could be replaced by the joint $\ell_s$-convergence

$$
(A_1^{(n)}, \ldots, A_m^{(n)}) \xrightarrow{\ell_s} (A_1^*, \ldots, A_m^*, b^*).
$$

REMARK 6.2.  We assume $s \le 2$ in the theorem above, unlike the finite-dimensional theorem in [27], where $s \le 3$. The reason is that for $2 < s \le 3$, we need to normalize the second moments, which is not generally possible in an infinite-dimensional Hilbert space.

**7. Contraction method for analytic functions.**  In this section, we link the general contraction theorem of Section 6 to recurrences of random analytic functions which will be needed subsequently. Suppose that $(X_n)_{n \ge 0}$ is a sequence of random analytic functions in a domain $D \subseteq \mathbb{C}$ such that, for some $n_0 \ge 1$,

$$
(7.1) \qquad X_n \overset{d}{=} \sum_{r=1}^{m} A_r^{(n)} \cdot X_{I_r^{(n)}}^{(r)} + b^{(n)}, \qquad n \ge n_0,
$$

where $A_1^{(n)}, \ldots, A_m^{(n)}$ and $b^{(n)}$ are random analytic functions in $D$ and $I^{(n)} = (I_1^{(n)}, \ldots, I_m^{(n)})$ is a vector of random integers with $I_r^{(n)} \in \{0, \ldots, n\}$. Furthermore, $\mathcal{L}(X_j^{(r)}) = \mathcal{L}(X_j)$ for all $r$ and $j$, and we have that $(A_1^{(n)}, \ldots, A_m^{(n)},$ $b^{(n)}, I^{(n)}), (X_j^{(1)}), \ldots, (X_j^{(m)})$ are independent. Note that the dot in (7.1) is meant as a pointwise complex multiplication. Hence, (7.1) is a special case of recurrence (6.3) with the random variables being analytic functions. The maps corresponding to $T$ in (6.1) now have, for some domain $\tilde{D}$, the form

$$(7.2) \qquad T : \mathcal{P}(\mathcal{H}(\tilde{D})) \to \mathcal{P}(\mathcal{H}(\tilde{D})),$$

$$\mathcal{L}(Z) \mapsto \mathcal{L}\left( \sum_{r=1}^m A_r^* \cdot Z^{(r)} + b^* \right),$$

where $A_1^*, \ldots, A_m^*$, and $b^*$ are random analytic functions in $\tilde{D}$, $\mathcal{L}(Z^{(r)}) = \mathcal{L}(Z)$ for $r = 1, \ldots, m$, and $(A_1^*, \ldots, A_m^*, b^*)$, $Z^{(1)}, \ldots, Z^{(m)}$ are independent.

We say that a random function $\Xi(z)$ is *locally bounded in $L^s$* if the function $z \mapsto \mathbb{E}|\Xi(z)|^s$ is locally bounded.

Theorem 6.1 then implies the following theorem.

Theorem 7.1. *Let $0 < s \le 2$. Let $(X_n)$ be as in (7.1) with $X_n$, $A_r^{(n)}$, $b^{(n)}$ being analytic functions of $z \in D$ locally bounded in $L^s$. For $1 < s \le 2$, assume that $\mathbb{E}X_n(z) = 0$ for all $n \ge 0$ and $z \in D$. Assume that $A_1^*, \ldots, A_m^*$ and $b^*$ are random analytic functions in $D$ and that $\Delta \subseteq D$ is a connected subset such that for each $x \in \Delta$, there exists a neighborhood $U_x \subseteq D$ of $x$ and a number $s(x) \le s$ such that*

$$(7.3) \qquad \sup_{z \in U_x} \mathbb{E}|A_r^{(n)}(z) - A_r^*(z)|^{s(x)} \to 0,$$

$$\sup_{z \in U_x} \mathbb{E}|b^{(n)}(z) - b^*(z)|^{s(x)} \to 0,$$

$$(7.4) \qquad \sup_{z \in U_x} \mathbb{E}|A_r^*(z)|^{s(x)} < \infty,$$

$$(7.5) \qquad \mathbb{E}\sum_{r=1}^m |A_r^*(x)|^{s(x)} < 1,$$

$$(7.6) \qquad \sup_{z \in U_x} \mathbb{E}[\mathbf{1}_{\{I_r^{(n)} \le \ell\}} |A_r^{(n)}(z)|^{s(x)}] \to 0,$$

*for all $\ell \in \mathbb{N}$ and $r = 1, \ldots, m$. There then exists a domain $\tilde{D} \subseteq D$ with $\Delta \subseteq \tilde{D}$ and a random analytic function $X$ in $\tilde{D}$ such that, as $n \to \infty$,*

$$(7.7) \qquad X_n \xrightarrow{d} X \qquad in \ \mathcal{H}(\tilde{D}),$$

*where $\mathcal{L}(X)$ is a fixed point of the map $T$ from (7.2). Moreover, for each $x \in \Delta$ there exists a neighborhood $\tilde{U}_x$ such that $\mathbb{E}\sup_{z \in \tilde{U}_x} |X(z)|^{s(x)} < \infty$ and, if $s(x) > 1$, also $\mathbb{E}X(z) = 0$, $z \in \tilde{U}_x$. $\mathcal{L}(X)$ is the only fixed point of $T$ with this property [even for a single $x \in \Delta$ and even with the integrability condition weakened to $\sup_{z \in \tilde{U}_x} \mathbb{E}|X(z)|^{s(x)} < \infty$].*

PROOF. We begin by showing that the assumptions hold locally uniformly in the following sense. Each $x \in \Delta$ has a neighborhood $U'_x = B(x, \varrho_x) \subseteq D$ with $\varrho_x > 0$ such that

$$\mathbb{E} \sup_{z \in U'_x} |A_r^{(n)}(z) - A_r^*(z)|^{s(x)} \to 0,$$

(7.8)

$$\mathbb{E} \sup_{z \in U'_x} |b^{(n)}(z) - b^*(z)|^{s(x)} \to 0,$$

(7.9)
$$\mathbb{E} \sum_{r=1}^m \sup_{z \in U'_x} |A_r^*(z)|^{s(x)} < 1,$$

(7.10)
$$\mathbb{E}\left[ \mathbf{1}_{\{I_r^{(n)} \le \ell\}} \sup_{z \in U'_x} |A_r^{(n)}(z)|^{s(x)} \right] \to 0.$$

To show this, we use the following simple lemma.

LEMMA 7.1. *If $\Psi$ is a random analytic function in a disc $B(w, r)$, then, for some constant $C > 0$ and every $s > 0$,*

$$\mathbb{E} \sup_{z \in B(w, r/2)} |\Psi(z)|^s \le C \sup_{z \in B(w, r)} \mathbb{E} |\Psi(z)|^s.$$

PROOF. Let $\gamma$ be the circle $|z - w| = 3r/4$, and let $|\gamma| = 3\pi r/2$ be its circumference. Since $|\Psi(z)|^s$ is subharmonic, standard properties of the Poisson kernel yield

$$\sup_{z \in B(w, r/2)} |\Psi(z)|^s \le C|\gamma|^{-1} \int_\gamma |\Psi(z)|^s \, dz.$$

Taking expectations on both sides implies the assertion. $\square$

Now, by assumption, for every $x \in \Delta$, there exists a disc $B(x, \varrho_x)$ such that (7.3) holds uniformly in $B(x, \varrho_x)$. Lemma 7.1 now shows that (7.8) holds with $U'_x = B(x, \varrho_x/2)$. Similarly, (7.6) and Lemma 7.1 applied to $\mathbf{1}_{\{I_r^{(n)} \le \ell\}} A_r^{(n)}$ yield (7.10). For (7.9), we first note that (7.4) and Lemma 7.1 imply that, for each $r = 1, \ldots, m$,

$$\mathbb{E} \sup_{z \in B(x, \varrho_x)} |A_r^*(z)|^{s(x)} < \infty$$

for a suitable $\varrho_x > 0$. Hence, dominated convergence shows that, as $\varrho \downarrow 0$, since $A_r^*$ is continuous,

$$\mathbb{E} \sup_{|z - x| \le \varrho} |A_r^*(z)|^{s(x)} \to \mathbb{E} |A_r^*(x)|^{s(x)}.$$

Summing over $r = 1, \ldots, m$ and using (7.5), we see that (7.9) holds for some $U'_x$.

Consider now a disc $U'_x$ such that (7.8)–(7.10) hold. We may assume that $\overline{U'_x} \subset D$. We let $H$ be the Bergman space $\mathcal{B}(U'_x)$ and regard (the restriction of) $X_n$ and $b^{(n)}$ as random elements of $H$, while the $A_r^{(n)}$ are regarded as random pointwise multiplication operators $f \mapsto A_r^{(n)} \cdot f$. Clearly, $\|A_r^{(n)}\|_{\mathrm{op}} \leq \sup_{z \in U'_x} |A_r^{(n)}(z)|$. Note that Lemma 7.1 implies that $\mathbb{E}\|X_n\|^s < \infty$, $\mathbb{E}\|b^{(n)}\|^s < \infty$, and $\mathbb{E}\|A_r^{(n)}\|_{\mathrm{op}}^s < \infty$.

By (7.8)–(7.10), Theorem 6.1 applies in $H$ [with $s$ replaced by $s(x)$] and shows that there exists a random analytic function $Z_x \in H = \mathcal{B}(U'_x)$ such that $X_n \overset{d}{\longrightarrow} Z_x$ in $\mathcal{B}(U'_x)$, and thus in $\mathcal{H}(U'_x)$, with $\mathcal{L}(Z_x)$ being the unique fixed point of $T$ from (7.2), with $\tilde{D} = U'_x$, in $\mathcal{P}_{s(x)}(\mathcal{B}(U'_x))$ for $0 < s(x) \leq 1$ and in $\mathcal{P}_{s(x),0}(\mathcal{B}(U'_x))$ for $1 < s(x) \leq 2$. The result now follows from Theorem 4.1 with $\tilde{D} = \bigcup_{x \in \Delta} U'_x$. Note that the condition $\mathbb{E}\|X\|_{\mathcal{B}(U'_x)}^{s(x)} = \mathbb{E}\|Z_x\|_{\mathcal{B}(U'_x)}^{s(x)} < \infty$, by Hölder's inequality, implies that

$$\mathbb{E}\int_{U'_x} |X(z)|^{s(x)} \, dm(z) = \mathbb{E}\|X\|_{L^{s(x)}(U'_x)}^{s(x)} \leq C(x)\mathbb{E}\|X\|_{L^2(U'_x)}^{s(x)} < \infty$$

and the fact that $|X(z)|^{s(x)}$ is subharmonic implies that for any strictly smaller disc $U''_x \subset U'_x$, we have $\sup_{z \in U''_x} |X(z)|^{s(x)} \leq C'(x) \int_{U'_x} |X(z)|^{s(x)} \, dm(z)$. Hence,

$$(7.11) \qquad \mathbb{E}\sup_{z \in U''_x} |X(z)|^{s(x)} < \infty$$

and thus

$$(7.12) \qquad \sup_{z \in U''_x} \mathbb{E}|X(z)|^{s(x)} < \infty.$$

Conversely, (7.11) implies that $\mathbb{E}\|X\|_{\mathcal{B}(U''_x)}^{s(x)} < \infty$ and (7.12) implies (7.11) in a smaller disc, by Lemma 7.1, so any of these, together with $\mathbb{E}X(z) = 0$ in a neighborhood of $x$ when $s(x) > 1$, yields uniqueness of $\mathcal{L}(X|_{\tilde{U}_x})$ for some neighborhood $\tilde{U}_x$, and thus uniqueness of $\mathcal{L}(X)$. $\quad\square$

REMARK 7.1. Theorem 7.1 and its proof immediately extend to (finite-dimensional) vector-valued functions $X_n$ and $b^{(n)}$ and matrix-valued functions $A_r^{(n)}$. Condition (7.5) then becomes $\mathbb{E}\sum_r \|A_r^*(x)\|_{\mathrm{op}}^s < 1$.

## 8. The expected profile.

The purpose of this section is to discuss in detail the first (and second) moment of the profile polynomials $W_n(z) = \sum_{k \geq 0} X_{n,k} z^k$ and the expected value $\mathbb{E}X_{n,k}$. We already know that for fixed $z \in D_1$ (recall the definition of the set $D_s$ in Section 3), the expected profile polynomial $\mathbb{E}W_n(z)$ behaves as $E(z)n^{\lambda_1(z)-1}$. In order to get more precise results, we need more information on $\lambda_1(z)$.

LEMMA 8.1.  *Suppose that $s > -t$. $D_s$ is then an open domain in the complex plane that contains the interval $(F(s), \infty)$ and $\lambda_1(z)$ is analytic in $D_s$.*

*Furthermore, if $z \in D_s$, then all $z' \in \mathbb{C}$ with $|z'| = |z|$ and $|\arg(z')| \leq |\arg(z)|$ are also contained in $D_s$. Moreover, the mapping $\varphi \mapsto \Re(\lambda_1(|z|e^{i\varphi}))$ is strictly decreasing for $0 \leq \varphi \leq |\arg(z)|$. In particular, $\Re(\lambda_1(z)) \leq \lambda_1(|z|)$, where equality holds if and only if $z$ is real and positive.*

PROOF.  Recall that $\lambda_1(z)$ is the root of $F(\lambda) = z$ with largest real part. Further, note that $F$ is a polynomial of degree $d = (m-1)(t+1)$. For the sake of brevity, we will only discuss the case $d > 2$. The three cases $m = 2$, $t = 0$ (where $d = 1$), $m = 3$, $t = 0$ (where $d = 2$) and $m = 2$, $t = 1$ (where $d = 2$) can be treated separately (and are, in fact, very easy).

We will first describe the "inverse map." For this purpose, we consider the mapping

$$\tau \mapsto F(\sigma + i\tau) \qquad (\tau \in \mathbb{R})$$

for fixed $\sigma > 0$. Since $F(\sigma - i\tau) = \overline{F(\sigma + i\tau)}$, it is sufficient to consider $\tau \geq 0$. By definition (3.5), it is clear that the argument $\arg(F(\sigma + i\tau))$ and the modulus $|F(\sigma + i\tau)|$ are strictly increasing functions for $\tau \geq 0$ and that we have $\lim_{\tau \to \infty} \arg(F(\sigma + i\tau)) = d\pi/2$. Hence, if $d > 2$, then there exists a unique minimal $\tau_0 = \tau_0(\sigma) > 0$ with $\arg(F(\sigma + i\tau_0)) = \pi$. Note that the mapping $\sigma \mapsto |F(\sigma + i\tau)|$ is also strictly increasing, but the mapping $\sigma \mapsto \arg(F(\sigma + i\tau))$ is strictly decreasing, for $\sigma > 0$ and fixed $\tau > 0$. This implies that the curves $\gamma_\sigma^+ := \{F(\sigma + i\tau) : 0 \leq \tau \leq \tau_0(\sigma)\}$, $\sigma > 0$, are all disjoint and that the mapping $\sigma \mapsto \tau_0(\sigma)$ is strictly increasing. Further, the mapping $\sigma \mapsto \tau_0(\sigma)$ is continuous. It also follows that we can parametrize $\gamma_\sigma^+$ as $\{r_\sigma(\phi)e^{i\phi} : 0 \leq \phi \leq \pi\}$ for some strictly increasing continuous function $r_\sigma$ on $[0, \pi]$.

The curve $\gamma_\sigma := \{F(\sigma + i\tau) : -\tau_0(\sigma) \leq \tau \leq \tau_0(\sigma)\} = \{r_\sigma(|\phi|)e^{i\phi} : -\pi \leq \phi \leq \pi\}$ is a simple closed curve that is the boundary of a compact set $K_\sigma = \{re^{i\phi} : r \in [0, r_\sigma(|\phi|)), \phi \in [-\pi, \pi]\}$, that is, the union of $\gamma_\sigma$ and its interior. Our next goal is to show that $D_s = \mathbb{C} \setminus (K_s \cup L_s)$, where $L_s := (-\infty, F(s + i\tau_0(s))]$ is a half-line. For this purpose, consider the set

$$Z_s := \{\lambda \in \mathbb{C} : \Re(\lambda) > s, -\tau_0(\Re(\lambda)) < \Im(\lambda) < \tau_0(\Re(\lambda))\}.$$

Suppose that $\lambda \in Z_s$. Then $F(\lambda) \in \gamma_{\Re(\lambda)}$ and thus $F(\lambda) \notin \gamma_s$. Moreover, $F(\lambda)$ can be connected to $\infty$ by a path disjoint from $\gamma_s$ (e.g., a piece of $\gamma_{\Re(\lambda)}$ plus the half-line $[F(\Re(\lambda)), \infty)$) and thus $F(\lambda)$ belongs to the exterior of $\gamma_s$, that is, $F(\lambda) \notin K_s$. Since $\lambda \in Z_s$ further implies that $|\arg(F(\lambda))| < \pi$, we also have $F(\lambda) \notin L_s$. Consequently, $F: Z_s \to \mathbb{C} \setminus (K_s \cup L_s)$.

Since the curves $\gamma_\sigma$, $\sigma > s$, are disjoint and simple, $F$ is injective on $Z_s$. Furthermore, to see that $F$ maps $Z_s$ onto $\mathbb{C} \setminus (K_s \cup L_s)$, suppose the contrary. Since $\mathbb{C} \setminus (K_s \cup L_s)$ is connected, there would be some $z \in \mathbb{C} \setminus (K_s \cup L_s)$ such that $z \in \overline{F(Z_s)} \setminus F(Z_s)$. Thus, there would exist a sequence $\lambda_n \in Z_s$ such that $F(\lambda_n) \to z$. This implies that the sequence $(\lambda_n)$ is bounded and there thus exists a subsequence converging to some $\lambda \in \overline{Z_s}$. By continuity, $F(\lambda) = z$, and since $z \notin F(Z_s)$, this implies that $\lambda \notin Z_s$ and thus that $\lambda \in \partial Z_s$. But, if $\lambda \in \partial Z_s$, then either $F(\lambda) \in \gamma_s \subset K_s$ or $\Im(\lambda) = \pm \tau_0(\Re(\lambda))$ and $F(\lambda) \in L_s$. Both cases contradict $F(\lambda) = z \in \mathbb{C} \setminus (K_s \cup L_s)$. Consequently, $F: Z_s \to \mathbb{C} \setminus (K_s \cup L_s)$ is a bijection.

Now, let $z \in \mathbb{C} \setminus (K_s \cup L_s)$. We will show that $z \in D_s$. We have just shown that then there exists $\theta \in Z_s$ with $F(\theta) = z$. By symmetry, we can assume that $\Im(\theta) \geq 0$. By the monotonicity properties of $|F(\sigma + i\tau)|$, it follows that $|F(\sigma + i\tau)| > |z|$ if $\sigma \geq \Re(\theta)$ and $|\tau| > \Im(\theta)$. Further, if $\sigma \geq \Re(\theta)$ and $|\tau| \leq \Im(\theta)$, then $\sigma + i\tau \in Z_s$, so $F(\sigma + i\tau) \neq F(\theta) = z$ unless $\sigma + i\tau = \theta$. Hence, $F(\lambda) = z$ has no other root with $\Re(\lambda) \geq \Re(\theta)$. Moreover, $F'(\theta) \neq 0$ (e.g., because $F$ is a bijection on $Z_s$) and thus $\theta$ is a simple root of $F(\lambda) = z$. Consequently, $\theta = \lambda_1(z)$ and $\Re(\lambda_1(z)) > \Re(\lambda_2(z))$, which implies that $z \in D_s$.

Similarly, if $\theta = s + i\tau'$ with $|\tau'| \leq \tau_0(s)$, then $|F(\sigma + i\tau)| > F(\theta)$ if $\sigma > s$ and $|\tau| \geq |\tau'|$, while if $\sigma > s$ and $|\tau| < |\tau'|$, then $|\arg(F(\sigma+i\tau))| < |\arg(F(\theta))|$. Hence, if $\sigma > s$ and $\tau \in \mathbb{R}$, then $F(\sigma + i\tau) \notin \gamma_s$. Since the half-plane $\{\sigma + i\tau : \sigma > s\}$ is connected, it is thus mapped by $F$ into the exterior of $\gamma_s$, that is, into $\mathbb{C} \setminus K_s$. Consequently, if $z \in K_s$, then $F(\lambda) = z$ has no root with $\Re\lambda > s$ and thus $z \notin D_s$. Finally, if $z \in L_s$, then $z = F(\sigma \pm i\tau_0(\sigma))$ for some $\sigma \geq s$. Thus, the fact that $z \in \gamma_\sigma$ and the argument just given together imply that $F(\theta) = z$ has no root with $\Re(\theta) > \sigma$. Hence, $\Re(\lambda_1(z)) = \sigma$, but there are two such roots, $\sigma \pm i\tau_0(\sigma)$, so $\Re(\lambda_1(z)) = \Re(\lambda_2(z))$ and $z \notin D_s$.

We have shown that

$$(8.1) \qquad D_s = \mathbb{C} \setminus (K_s \cup L_s) = \{re^{i\phi} : r > r_s(|\phi|), -\pi < \phi < \pi\}$$

and that the inverse mapping $F^{-1}: \mathbb{C} \setminus (K_s \cup L_s) \to Z_s$ explicitly computes $\lambda_1(z) = F^{-1}(z)$, which is a simple root. Note, too, that $\lambda_1(z)$ (for $z \in D_s$) is characterized by the property that it has smallest absolute imaginary part among all solutions of $F(\lambda) = z$ with $\Re(\lambda) > s$. By the implicit function theorem, $\lambda_1(z)$ is analytic in $D_s$.

Since $F(s)$ is the only boundary point of $D_s$ on the positive real line, it follows that $D_s$ contains the interval $(F(s), \infty)$. [Alternatively, use (8.1).]

Finally, (8.1) and the fact that $r_s$ is strictly increasing on $[0, \pi]$ together imply that if $re^{i\phi} \in D_s$ and $|\phi'| < |\phi| < \pi$, then $r > r_s(|\phi|) > r_s(|\phi'|)$ and $re^{i\phi'} \in D_s$. Moreover, if $\sigma := \Re(\lambda_1(re^{i\phi}))$, then $re^{i\phi} \in \partial D_\sigma = \gamma_\sigma \cup L_\sigma$, and thus $re^{i\phi} \in \gamma_\sigma$. The same argument shows that $r = r_\sigma(|\phi|) > r_\sigma(|\phi'|)$ and $re^{i\phi'} \in D_\sigma$. Hence, $\Re(\lambda_1(re^{i\phi'})) > \sigma = \Re(\lambda_1(re^{i\phi}))$ and the final statement in the lemma follows. $\square$

The next step is an extension of Lemma 3.1. Note that $\Re(\lambda_1(z))$ is well defined for all $z \in \mathbb{C}$. (Estimates involving $\log n$ are only supposed to hold for $n \geq 2$.)

LEMMA 8.2. *Let $W_n(z) = \sum_{k \geq 0} X_{n,k} z^k$ denote the (random) profile polynomials.*

(i) *If $K$ is a compact subset of $D_1$, then there exists $\delta > 0$ such that*

$$(8.2) \qquad \mathbb{E}W_n(z) = n^{\lambda_1(z)-1}(E(z) + O(n^{-\delta}))$$

*uniformly for $z \in K$.*

(ii) *if $K$ is a compact subset of $\mathbb{C}$, then there exists $D \geq 0$ such that*

$$(8.3) \qquad |\mathbb{E}W_n(z)| \lesssim n^{\max\{\Re(\lambda_1(z))-1,0\}}(\log n)^D$$

*uniformly for $z \in K$.*

PROOF. The proof is a direct extension of the results of [9] applied to the recurrence relation (3.3). In particular, we have to take care of the uniformity in $z \in K$. This can be done by a careful inspection of the proof in [9]; see the Appendix. $\square$

With the help of Lemma 8.2, we directly obtain bivariate asymptotic expansions for $\mathbb{E}X_{n,k}$ in a large range. It turns out that one must solve the equation

$$(8.4) \qquad \beta\lambda_1'(\beta) = \alpha.$$

From (3.5), it follows that

$$(8.5) \qquad \begin{aligned} \beta\lambda_1'(\beta) &= F(\lambda_1(\beta))/F'(\lambda_1(\beta)) \\ &= \left( \frac{1}{\lambda_1(\beta)+t} + \frac{1}{\lambda_1(\beta)+t+1} + \cdots + \frac{1}{\lambda_1(\beta)+(t+1)m-2} \right)^{-1}. \end{aligned}$$

Note that this formula also shows that the mapping $\beta \mapsto \beta\lambda_1'(\beta)$ is strictly increasing because $\lambda_1(\beta)$ is strictly increasing for $\beta > 0$. Moreover, $\lambda_1(\beta)$ increases from $-t$ to $\infty$ for $0 < \beta < \infty$ and it follows that $\beta\lambda_1'(\beta)$ increases from 0 to $\infty$. Hence, (8.4) has a unique solution $\beta(\alpha) > 0$ for every $\alpha > 0$, with $\beta(\alpha)$ strictly increasing. Since we must assume that $\lambda_1 > 1$, we note that (8.4) has a proper solution with $\lambda_1(\beta) > 1$, and thus $\beta \in D_1$, if and only if

$$\alpha > \alpha_0 := \left( \frac{1}{t+1} + \frac{1}{t+2} + \cdots + \frac{1}{(t+1)m-1} \right)^{-1}.$$

LEMMA 8.3.  *Suppose that* $\alpha_1, \alpha_2$ *with* $\alpha_0 < \alpha_1 < \alpha_2 < \infty$ *are given and let* $\beta(\alpha)$ *be defined by* $\beta(\alpha)\lambda_1'(\beta(\alpha)) = \alpha$. *Then*

$$\mathbb{E} X_{n,k} = \frac{E(\beta(\alpha_{n,k}))n^{\lambda_1(\beta(\alpha_{n,k})) - \alpha_{n,k} \log(\beta(\alpha_{n,k})) - 1}}{\sqrt{2\pi(\alpha_{n,k} + \beta(\alpha_{n,k})^2 \lambda_1''(\beta(\alpha_{n,k}))) \log n}} (1 + O((\log n)^{-1/2}))$$

*uniformly for* $\alpha_{n,k} = k/\log n \in [\alpha_1, \alpha_2]$ *as* $n, k \to \infty$.

PROOF.  By Cauchy's formula, we have

$$\mathbb{E} X_{n,k} = \frac{1}{2\pi \mathrm{i}} \int_{|z|=\beta} \mathbb{E} W_n(z) z^{-k-1} \, dz.$$

Note that $\mathbb{E} W_n(z) z^{-k}$ behaves as

$$(8.6) \qquad \mathbb{E} W_n(z) z^{-k} \sim E(z) n^{\lambda_1(z)-1} z^{-k} = \frac{1}{n} E(z) e^{\lambda_1(z) \log n - k \log z}.$$

In order to evaluate the above Cauchy integral, we use a standard saddle point method. The saddle point of the function $z \mapsto \lambda_1(z) \log n - k \log z$ is given by $z_0 = \beta$, which satisfies $\beta \lambda_1'(\beta) = k/\log n$, that is, by $\beta(k/\log n)$.

By construction, the real interval $[\beta(\alpha_1), \beta(\alpha_2)]$ is contained in $D_1$. Hence, there exists $\gamma > 0$ such that the set $\{z \in \mathbb{C} : |z| \in [\beta(\alpha_1), \beta(\alpha_2)], |\arg(z)| \le \gamma\}$ is also contained in $D_1$.

Let $K = \{z \in \mathbb{C} : |z| \in [\beta(\alpha_1), \beta(\alpha_2)], \gamma \le |\arg(z)| \le \pi\}$. By Lemma 8.1, there then exists $\eta > 0$ such that for all $\beta \in [\beta(\alpha_1), \beta(\alpha_2)]$,

$$\max_{z \in K, |z|=\beta} \max\{\Re(\lambda_1(z)), 1\} \le \lambda_1(\beta) - \eta.$$

[Uniformity follows from the continuity of $\lambda_1(z)$.] Hence, by Lemma 8.2,

$$\int_{|z|=\beta, |\arg(z)| \ge \gamma} |\mathbb{E} W_n(z) z^{-k-1}| \, dz \lesssim n^{\lambda_1(\beta(\alpha)) - 1 - \eta/2 - \alpha \log(\beta(\alpha))},$$

where $\alpha = k/\log n$ and $\beta = \beta(\alpha)$. Thus, this part of the integral is negligible.

For the remaining integral (leading to the asymptotic leading term), we use the substitution $z = \beta e^{\mathrm{i}t}$ ($|t| \le \gamma$) and the approximation

$$(8.7) \qquad \lambda_1(z) \log n - k \log z = (\lambda_1(\beta) - \alpha \log \beta) \log n$$

$$(8.8) \qquad\qquad\qquad + \frac{1}{2}((\lambda_1''(\beta) + \alpha\beta^{-2}) \log n)(z - \beta)^2$$

$$(8.9) \qquad\qquad\qquad + O(\log n |z - \beta|^3)$$

to obtain the final form after standard *saddle point algebra.*  □

In what follows, we will also need estimates for the second moments of $W_n(z)$.

LEMMA 8.4.    *For every compact set $K \subseteq \mathbb{C}$, we have*

$$(8.10) \qquad \mathbb{E}|W_n(z)|^2 = O(n^{\max\{\lambda_1(|z|^2)-1, 2\Re(\lambda_1(z))-2, 0\}} (\log n)^{D'})$$

*uniformly for $z \in K$, where $D' \geq 0$ is an absolute constant.*

PROOF.    We use (3.2) twice, for $z$ and $\overline{z}$, and obtain

$$
\begin{aligned}
|W_n(z)|^2 &\stackrel{\mathrm{d}}{=} |z|^2 (|W_{V_{n,1}}^{(1)}(z)|^2 + \cdots + |W_{V_{n,m}}^{(m)}(z)|^2) \\
(8.11) &\qquad + |z|^2 \sum_{i \neq j} W_{V_{n,i}}^{(i)}(z) W_{V_{n,j}}^{(j)}(\overline{z}) \\
&\qquad + 2(m-1)\Re\left(\sum_{j=1}^m z W_{V_{n,j}}^{(j)}(z)\right) + (m-1)^2.
\end{aligned}
$$

We take the expectation. By Lemma 8.2, $|\mathbb{E}W_n(z)| = O(A_n(z))$, uniformly for $z \in K$, where $A_n(z) := n^{\max\{\Re(\lambda_1(z))-1, 0\}} (\log n)^D$ (for some fixed $D \geq 0$) for $n \geq 2$, and $A_0(z) := A_2(z) := 1$, say. Hence, for $l < n$, and uniformly in $z \in K$,

$$\mathbb{E}(W_{V_{n,1}}(z) \mid V_{n,1} = l) = O(A_l(z)) = O(A_n(z))$$

and thus

$$\mathbb{E}(zW_{V_{n,1}}(z) \mid V_{n,1} = l) = O|\mathbb{E}(W_{V_{n,1}}(z) \mid V_{n,1} = l)| = O(A_n(z)) = O(A_n(z)^2).$$

Similarly, for $l_1 + l_2 < n$, and uniformly in $z \in K$,

$$\mathbb{E}(W_{V_{n,1}}(z) W_{V_{n,2}}(z) \mid V_{n,1} = l_1, V_{n,2} = l_2) = O(A_{l_1}(z) A_{l_2}(z)) = O(A_n(z)^2).$$

Consequently, uniformly for $x \in K$, (8.11) yields, using (2.2),

$$(8.12) \qquad \mathbb{E}|W_n(z)|^2 = m|z|^2 \mathbb{E}|W_{V_{n,1}}(z)|^2 + O(A_n(z)^2)$$

$$(8.13) \qquad = m|z|^2 \sum_{\ell=0}^{n-1} \frac{\binom{\ell}{t}\binom{n-\ell-1}{(m-1)t+m-2}}{\binom{n}{mt+m-1}} \mathbb{E}|W_\ell(z)|^2$$

$$(8.14) \qquad + O(n^{\max\{2\Re(\lambda_1(z))-2, 0\}} (\log n)^{2D}).$$

This is an equation of the same type as (3.3) and we can again apply [9] to obtain the stated estimate. As in the proof of Lemma 8.2, an inspection of [9] shows that the estimate holds uniformly in $z$; see the Appendix for details.    □

REMARK 8.1. A special case of this result for $m = 2$ and $t = 0$ has been proven in [4]. In this case, we have $\lambda_1(z) = 2z$ and obtain (for some $D \geq 0$; in fact, $D = 2$ will suffice for all $z$ and $D = 1$ or $0$ will suffice for all $z \neq 1/2$)

$$\mathbb{E}|W_n(z)|^2 \lesssim n^{\max\{4\Re z - 2, 0\}} (\log n)^D \qquad (|z - 1| \leq 1/\sqrt{2})$$

and

$$\mathbb{E}|W_n(z)|^2 \lesssim n^{\max\{2|z|^2 - 1, 0\}} (\log n)^D \qquad (|z - 1| \geq 1/\sqrt{2}).$$

REMARK 8.2. The method of Lemma 8.4 can be used for many other *functionals* of $W_n(z)$. For example, the expected derivative $\mathbb{E}W_n'(z)$ satisfies the recurrence

$$\mathbb{E}W_n'(z) = mz \sum_{\ell=0}^{n-1} \frac{\binom{\ell}{t}\binom{n-\ell-1}{(m-1)t+m-2}}{\binom{n}{mt+m-1}} \mathbb{E}W_\ell'(z)$$

$$+ m \sum_{\ell=0}^{n-1} \frac{\binom{\ell}{t}\binom{n-\ell-1}{(m-1)t+m-2}}{\binom{n}{mt+m-1}} \mathbb{E}W_\ell(z).$$

For simplicity, let $z \in D_1$ be real and nonnegative. Then, from

$$m \sum_{\ell=0}^{n-1} \frac{\binom{\ell}{t}\binom{n-\ell-1}{(m-1)t+m-2}}{\binom{n}{mt+m-1}} \mathbb{E}W_\ell(z) = \frac{1}{z}(\mathbb{E}W_n(z) - (m-1)) = O(n^{\lambda_1(z)-1})$$

and an application of [9], we obtain

$$\mathbb{E}W_n'(z) = O(n^{\lambda_1(z)-1} \log n).$$

[This also follows from (8.2) by Cauchy's estimates.]

We close this section with a proof that the sets $I$ and $I'$ (defined in Theorem 1.1) are in fact intervals.

LEMMA 8.5. Let $I := \{\beta > 0 : 1 < \lambda_1(\beta^2) < 2\lambda_1(\beta) - 1\}$ and $I' := \{\beta\lambda_1'(\beta) : \beta \in I\}$. Then $I$ and $I'$ are open intervals that are contained in the positive real line. More precisely, $1 \in I \subseteq (\frac{1}{m}, \beta(\alpha_+))$ and $\alpha_{\max} \in I' \subseteq (\alpha_0, \alpha_+)$.

PROOF. Since $\lambda_1(z)$ is increasing for $z > 0$, it is clear that $I_1 := \{\beta > 0 : 1 < \lambda_1(\beta^2)\}$ is an interval. We show that $I_2 := \{\beta > 0 : \lambda_1(\beta^2) < 2\lambda_1(\beta) - 1\}$ is also an interval, which implies that $I = I_1 \cap I_2$ is an interval.

Suppose that $\beta > 0$ and that $\lambda = \lambda_1(\beta)$ with $\lambda_1(\beta^2) = 2\lambda_1(\beta) - 1$. Then $F(2\lambda - 1) = F(\lambda_1(\beta^2)) = \beta^2 = F(\lambda)^2$. However, for $\lambda > -(t-1)/2$, we have

$$\frac{d}{d\lambda}(\log F(2\lambda - 1) - 2\log F(\lambda)) = \sum_{i=t}^{mt+m-2} \left(\frac{2}{2\lambda - 1 + i} - \frac{2}{\lambda + i}\right),$$

which is greater than 0 for $\lambda < 1$ and less than 0 for $\lambda > 1$. Thus, $q(\lambda) := F(2\lambda - 1)/F(\lambda)^2$ is strictly increasing on $[(1-t)/2, 1]$ and strictly decreasing on $[1, \infty)$. Moreover, $q((1-t)/2) = 0$, $q(1) = 1/F(1) = m > 1$ and $q(\lambda) \to 0$ as $\lambda \to \infty$. Consequently, there are exactly two roots, $\lambda_1^* < \lambda_2^*$, of $F(2\lambda - 1) = F(\lambda)^2$ in $[(1-t)/2, \infty)$ and two roots $\beta_j^* = F(\lambda_j^*) > 0$ of $\lambda_1(\beta^2) = 2\lambda_1(\beta) - 1$. [Note that $2\lambda_1(\beta) - 1 = \lambda_1(\beta^2) > -t$ implies that $\lambda_1(\beta) > (1-t)/2$.] Since $\lambda_1(1) = 2$, it is easily seen that $\lambda_1(\beta^2) < 2\lambda_1(\beta) - 1$ (for $\beta > 0$) if and only if $\beta_1^* < \beta < \beta_2^*$.

Set $I = (\underline{\beta}, \overline{\beta})$. Since $\lambda(1) = 2$, we surely have $1 \in I$. Next, note that $\beta = \frac{1}{m}$ corresponds to $\lambda_1(\beta) = 1$. Thus, $\lambda_1(1/m^2) < \lambda_1(1/m) = 1$, which implies that $\frac{1}{m} < \underline{\beta}$.

In order to prove $\beta_+ := \beta(\alpha_+) > \overline{\beta}$, it suffices to show that $\lambda_1(\beta_+^2) > 2\lambda_1(\beta_+) - 1$ or, equivalently, that $F(2\lambda_1(\beta_+) - 1) < F(\lambda_1(\beta_+))^2$ [since $F(\lambda_1(\beta_+))^2 = \beta_+^2 = F(\lambda_1(\beta_+^2))$]. First, by definition,

$$\log F(\lambda_1(\beta_+)) = (\lambda_1(\beta_+) - 1) \sum_{i=t}^{mt+m-2} \frac{1}{\lambda_1(\beta_+) + i}.$$

Moreover, with $S_+ := \sum_{i=t}^{mt+m-2} (\lambda_1(\beta_+) + i)^{-1}$, it follows by a convexity argument (compare with [5], Lemma 3.2) that, for every $\lambda \geq 1$,

$$\log F(\lambda) \leq (\lambda - 1) S_+,$$

with equality only for $\lambda = \lambda_1(\beta_+)$. Consequently,

$$\log F(2\lambda_1(\beta_+) - 1) < (2\lambda_1(\beta_+) - 2) S_+ = 2 \log F(\lambda_1(\beta_+)).$$

Thus, we have $F(2\lambda_1(\beta_+) - 1) < F(\lambda_1(\beta_+))^2$ and consequently $\beta(\alpha_+) > \overline{\beta}$.

Finally, since the mapping $\beta \mapsto \beta \lambda_1'(\beta)$ is strictly increasing, by (8.5), and continuous, it also follows that $I'$ is an interval [that is contained in $(\alpha_0, \alpha_+)$]. $\quad\square$

REMARK 8.3. In Theorem 11.1, we consider the set $J = \{\beta > 0 : \lambda_1(\beta^2) < 2\lambda_1(\beta) - 1\}$ instead of $I$. This equals $I_2$ in the proof above, so $J$ is also an open interval. Furthermore, a slight extension of the above proof shows that $J \subset (\beta(\alpha_-), \beta(\alpha_+))$. The proof shows that $\beta_+ > \sup J$, and $\beta_- := \beta(\alpha_-) < \inf J$ can be shown in exactly the same way.

**9. Proof of Theorem 1.2.** In this section, we prove Theorem 1.2.

PROOF OF THEOREM 1.2. The sequence of random analytic functions $(W_n)$ in Theorem 1.2 satisfies the recurrence (3.2). Hence, for

$$X_n(z) := \frac{W_n(z) - \mathbb{E}W_n(z)}{\mathbb{E}W_n(z)} = \frac{W_n(z)}{\mathbb{E}W_n(z)} - 1,$$

we obtain, with $G_n(z) := \mathbb{E}W_n(z)$, that

$$X_n(z) \stackrel{d}{=} \sum_{r=1}^{m} z \frac{G_{V_{n,r}}(z)}{G_n(z)} X_{V_{n,r}}^{(r)} + \frac{1}{G_n(z)} \left( m - 1 - G_n(z) + z \sum_{r=1}^{m} G_{V_{n,r}}(z) \right).$$

Hence, we have

$$X_n \stackrel{d}{=} \sum_{r=1}^{m} A_r^{(n)} \cdot X_{I_r^{(n)}}^{(r)} + b^{(n)},$$

with $I_r^{(n)} = V_{n,r}$, $A_r^{(n)} = zG_{V_{n,r}}(z)/G_n(z)$,

$$b^{(n)} = \frac{1}{G_n(z)} \left( m - 1 - G_n(z) + z \sum_{r=1}^{m} G_{V_{n,r}}(z) \right)$$

and conditions as in (7.1). We will see below that the sequence $(X_n)$ of random analytic functions satisfies the conditions of Theorem 7.1 for all $1 < s \le 2$ with $D = \{z \in D_1 : E(z) \ne 0\}$,

$$A_r^*(z) = zV_r^{\lambda_1(z)-1}, \qquad b^* = z \sum_{r=1}^{m} V_r^{\lambda_1(z)-1} - 1,$$

for $r = 1, \ldots, m$, and $\Delta = (1/m, \beta(\alpha_+))$.

Theorem 7.1 then implies that $X_n \stackrel{d}{\longrightarrow} X$ in $\mathcal{H}(\tilde{D})$, where $\tilde{D}$ is a complex neighborhood of the real interval $(1/m, \beta(\alpha_+))$ and $\mathcal{L}(X)$ is the fixed point of $T$ defined in (7.2), with the integrability condition from Theorem 7.1. Recall that for $x \in (1/m, \beta(\alpha_+))$, we have $\lambda_1(x) > 1$ and note that this convergence implies the assertion since $W_n(z)/\mathbb{E}W_n(z) = X_n + 1$. Hence, we have

$$\frac{W_n(z)}{\mathbb{E}W_n(z)} \stackrel{d}{\longrightarrow} Y(z) = X(z) + 1 \qquad \text{in } \mathcal{H}(\tilde{D}),$$

where

$$Y \stackrel{d}{=} \sum_{r=1}^{m} zV_r^{\lambda_1(z)-1} \cdot Y^{(r)},$$

with conditions as in (7.2), which is (3.7). The integrability condition on $X$ is obviously equivalent to the same condition for $Y$ and since, as we shall see below, we may take $s(x) > 1$ arbitrarily close to 1, the condition is equivalent to the existence, for each $x \in I$, of some $s(x) > 1$ such that $\mathbb{E}|Y(z)|^{s(x)}$ is finite and bounded in a neighborhood of $x$, as asserted in Section 1.

It remains to verify conditions (7.3)–(7.6). Using Lemma 8.2, we obtain, uniformly in each compact subset of $D$,

$$A_r^{(n)} z \frac{G_{V_{n,r}}(z)}{G_n(z)} = z \frac{V_{n,r}^{\lambda_1(z)-1}(E(z) + O(V_{n,r}^{-\delta}))}{n^{\lambda_1(z)-1}(E(z) + O(n^{-\delta}))}$$

$$= z \left( \frac{V_{n,r}}{n} \right)^{\lambda_1(z)-1} (1 + O(V_{n,r}^{-\delta})).$$

[$V_{n,r}$ may equal 0, but that is not a problem; the cautious reader may write $(1 + V_{n,r})^{-\delta}$ above.]

We have $V_{n,r}/n \xrightarrow{\text{d}} V_r$ by (3.6) and thus we may, by a suitable coupling, assume that $V_{n,r}/n \to V_r$ a.s. (see also Remark 3.1). Since these random variables are bounded by 1 and $\Re\lambda_1(z) - 1 > 0$ in $D$, dominated convergence yields $A_r^{(n)}(z) \to A_r^*(z)$ in $L^s$ for any $s > 0$. This also implies that $b^{(n)}(z) \to b^*(z)$ in $L^s$. Moreover, these $L^s$-convergences are uniform in any compact subset of $D$ and arbitrary $s > 1$. This establishes condition (7.3). For bounded neighborhoods $U_x$ and arbitrary $1 < s(x) \le 2$, we have that $|A_r^*(z)|^{s(x)}$ is uniformly bounded in $z \in U_x$. This implies conditions (7.4) and (7.6), since we have $\mathbb{P}(I_r^{(n)} \le \ell) \to 0$ for all $\ell \in \mathbb{N}$ as $I_r^{(n)}/n \to V_r$ and $\mathbb{P}(V_r = 0) = 0$.

For condition (7.5), note that $V_r$ has the $\mathrm{Beta}(t+1, (m-1)(t+1))$ distribution. This implies that, for $\alpha > 0$, we have

$$\mathbb{E}V_r^\alpha = \frac{\Gamma(t + \alpha + 1)((t+1)m-1)!}{t!\Gamma(m(t+1) + \alpha)} = \frac{1}{mF(\alpha+1)},$$

with $F$ given in (3.5). Let $x \in \Delta = (1/m, \beta(\alpha_+))$. We have

$$\sum_{r=1}^m \mathbb{E}|A_j^*(x)|^s = \frac{x^s}{F(s\lambda_1(x) - s + 1)} =: g_x(s).$$

We have $g_x(1) = 1$. Thus, the existence of an $s(x) \in (1, s)$ with (7.5) follows from $g_x'(1) < 0$. To verify a negative derivative, we consider $h_x(s) := -\log(g_x(s)) = \log(F(s(\lambda_1(x) - 1) + 1)) - s\log x$. We then have

$$h_x'(1) = (\lambda_1(x) - 1)\frac{F'}{F}(\lambda_1(x)) - \log x$$

$$= \sum_{i=t}^{mt+m-2} \frac{\lambda_1(x) - 1}{\lambda_1(x) + i} - \log x.$$

From $\log x = \log(F(\lambda_1(x)))$ and (1.4), we obtain that the only zeros of $h_x'(1)$ are at $x = F(\lambda_-) = \beta(\alpha_-)$ and $x = F(\lambda_+) = \beta(\alpha_+)$. For $x = 1$, we obtain, with $\lambda_1(1) = 2$, that $h_1'(1) > 0$, thus, by continuity of $x \mapsto h_x'(1)$, we obtain $h_x'(1) > 0$ for all $\beta(\alpha_-) < x < \beta(\alpha_+)$.

Thus, for all $\beta(\alpha_-) < x < \beta(\alpha_+)$, there exists $s(x) \in (1, s)$ such that $g_x(s(x)) < 1$. In particular, this shows (7.5). We have verified the conditions of Theorem 7.1 and the proof is complete.  $\square$

## 10. Reduction to the profile.

We now come back to the original problem. We know (Theorem 1.2) that the profile polynomials $W_n(z)$ satisfy a functional limit theorem

$$(10.1) \qquad (W_n(z)/\mathbb{E}W_n(z), z \in B) \xrightarrow{\text{d}} (Y(z), z \in B)$$

for some open domain $B \subseteq \mathbb{C}$ including the open interval $(1/m, \beta(\alpha_+)) \subseteq \mathbb{R}$ and that $(Y(z), z \in B)$ is the process (of random analytic functions) that satisfies the stochastic fixed point equation (3.7) with $\mathbb{E}Y(z) = 1$ and a certain integrability condition.

The idea is now to reconstruct $X_{n,k}$ from the limit relation (10.1). In particular, we want to show that

$$(10.2) \qquad \left( \frac{X_{n,\lfloor \alpha \log n \rfloor}}{\mathbb{E}X_{n,\lfloor \alpha \log n \rfloor}}, \alpha \in I' \right) \xrightarrow{\mathrm{d}} (Y(\beta(\alpha)), \alpha \in I'),$$

where $I'$ (defined in Theorem 1.1) has the property that all $\beta = \beta(\alpha)$ (for $\alpha \in I'$) satisfy $1 < \lambda_1(\beta^2) < 2\lambda_1(\beta) - 1$.

We use the Cauchy formula and split it into two parts. More precisely, we fix a compact interval $I_c \subseteq (1/m, \beta(\alpha_+))$ and a small $\varphi > 0$ such that the compact set $B_1 := \{z \in \mathbb{C} : |z| \in I_c, |\arg(z)| \leq \varphi\}$ is contained in $B$. Further, let $I_c' := \{\beta\lambda_1'(\beta) : \beta \in I_c\}$ and $B_2 := \{z \in \mathbb{C} : |z| \in I_c, \varphi < |\arg(z)| \leq \pi\}$.

For $\alpha \in I_c'$, so that $\beta(\alpha) \in I_c$, we write

$$(10.3) \qquad \begin{aligned} X_{n,\lfloor \alpha \log n \rfloor} &= \frac{1}{2\pi \mathrm{i}} \int_{|z| = \beta(\alpha), z \in B_1} W_n(z) z^{-\lfloor \alpha \log n \rfloor - 1} \, dz \\ &\quad + \frac{1}{2\pi \mathrm{i}} \int_{|z| = \beta(\alpha), z \in B_2} W_n(z) z^{-\lfloor \alpha \log n \rfloor - 1} \, dz. \end{aligned}$$

We study the two integrals separately. For the first part, we define linear operators $T_n$, mapping the space $C(B_1)$ of continuous functions on $B_1$ into the space $D(I_c')$ of right-continuous functions with left limits on $I_c'$, by

$$(10.4) \qquad T_n(G)(\alpha) = \frac{1/(2\pi i) \int_{|z| = \beta(\alpha), z \in B_1} G(z) \mathbb{E}W_n(z) z^{-\lfloor \alpha \log n \rfloor - 1} \, dz}{\mathbb{E}X_{n,\lfloor \alpha \log n \rfloor}},$$

$$\alpha \in I_c'.$$

Note that if we take $G(z) = W_n(z)/\mathbb{E}W_n(z)$, the numerator in (10.4) equals the first term on the right-hand side of (10.3). The second term will be shown to be small and thus $T_n(W_n(z)/\mathbb{E}W_n(z))$ is an approximation of

$$(X_{n,\lfloor \alpha \log n \rfloor}/\mathbb{E}X_{n,\lfloor \alpha \log n \rfloor}, \alpha \in I_c').$$

We begin by studying $T_n$ in Lemma 10.1. We will then, in Lemma 10.2, show that the second term of (10.3) is sufficiently small to be neglected and Theorem 1.1 will follow.

We will use the supremum norm; for convenience, we write, from any set $E$,

$$\|f\|_E := \sup_E |f|.$$

LEMMA 10.1. (i) *The operators $T_n$ are uniformly continuous with respect to the supremum norm. More precisely, there exists a constant $C > 0$ (depending on $I_c$ and $B_1$) such that*

$$\|T_n(F) - T_n(G)\|_{I'_c} \leq C \cdot \|F - G\|_{B_1}.$$

(ii) *If $F_n \to F$ uniformly on $B_1$, then $T_n(F_n) \to F$ uniformly on $I'_c$.*

PROOF. (i) Suppose that $\|F - G\|_{B_1} \leq \delta$ and that $\alpha \in I'_c$. Then

$$\begin{aligned}
(10.5) \quad & |T_n(F)(\alpha) - T_n(G)(\alpha)| \\
& \leq \frac{\delta}{\mathbb{E}X_{n,\lfloor \alpha \log n \rfloor}} \frac{1}{2\pi} \int_{|z|=\beta(\alpha), z \in B_1} |\mathbb{E}W_n(z) z^{-\lfloor \alpha \log n \rfloor - 1}||dz|.
\end{aligned}$$

First, suppose that $k = \alpha \log n$ is an integer. Lemma 8.3 yields an estimate of $\mathbb{E}X_{n,k}$ and its proof, in particular (8.6) and (8.7), yields an estimate of the same order for the integral in (10.5). Hence, (10.5) implies that

$$(10.6) \qquad |T_n(F)(\alpha) - T_n(G)(\alpha)| \leq C\delta$$

for some $C$, uniformly in $\alpha \in I'_c$, such that $\alpha \log n$ is an integer.

For general $\alpha$, we define $\alpha' = \lfloor \alpha \log n \rfloor / \log n$ and note that $|\alpha' - \alpha| \leq 1/\log n$, so $|\beta(\alpha') - \beta(\alpha)| = O(1/\log n)$. It is easily checked that if we replace $\alpha'$ by $\alpha$ in the estimate of $X_{n,k}$ in Lemma 8.3, then the result will change by at most a factor $n^{O(1/\log n)} = e^{O(1)}$. It follows that

$$(10.7) \qquad \mathbb{E}X_{n,\lfloor \alpha \log n \rfloor} \gtrsim (\log n)^{-1/2} n^{\lambda_1(\beta(\alpha)) - \alpha \log \beta(\alpha) - 1}$$

and that (10.6) holds uniformly in all $\alpha \in I'_c$, possibly with a larger constant.

For (ii), let $F \in C(B_1)$. By using standard saddle point techniques, as in the proof of Lemma 8.3, we have

$$\lim_{n \to \infty} \frac{1}{\mathbb{E}X_{n,\lfloor \alpha \log n \rfloor}} \frac{1}{2\pi i} \int_{|z|=\beta(\alpha), z \in B_1} F(z) \mathbb{E}W_n(z) z^{-\lfloor \alpha \log n \rfloor - 1} \, dz = F(\beta(\alpha)),$$

that is, $T_n(F) \to F$, uniformly on $I'_c$. Finally, if $F_n \to F$ uniformly on $B_1$, then $T_n(F_n) - T_n(F) \to 0$, by (i), and consequently $T_n(F_n) \to F$ uniformly on $I'_c$. This completes the proof of the lemma. $\square$

Next, we focus on the error

$$\begin{aligned}
(10.8) \quad & \frac{X_{n,\lfloor \alpha \log n \rfloor}}{\mathbb{E}X_{n,\lfloor \alpha \log n \rfloor}} - T_n \left( \frac{W_n(z)}{\mathbb{E}W_n(z)} \right) \\
& = \frac{1/(2\pi i) \int_{|z|=\beta(\alpha), z \in B_2} W_n(z) z^{-\lfloor \alpha \log n \rfloor - 1} \, dz}{\mathbb{E}X_{n,\lfloor \alpha \log n \rfloor}},
\end{aligned}$$

where we recall that $B_2 = \{z \in \mathbb{C} : |z| \in I_c, \varphi < |\arg(z)| \leq \pi\}$.

Lemma 10.2. *For every compact interval $I_c$ contained in $I = \{\beta > 0 : 1 < \lambda_1(\beta^2) < 2\lambda_1(\beta) - 1\}$,*

$$\sup_{\alpha \in I_c'} \left| \frac{1/(2\pi) \int_{|z|=\beta(\alpha), z \in B_2} W_n(z) z^{-\lfloor \alpha \log n \rfloor - 1} \, dz}{\mathbb{E} X_{n, \lfloor \alpha \log n \rfloor}} \right| \xrightarrow{\mathrm{p}} 0$$

*as $n \to \infty$.*

Proof. Let

$$G_{n,\alpha}(z) := \frac{|W_n(z) z^{-\lfloor \alpha \log n \rfloor}|}{\mathbb{E} X_{n, \lfloor \alpha \log n \rfloor}}.$$

If we further define

$$(10.9) \qquad H_n(z) := \frac{W_n(z)}{n^{\lambda_1(|z|)-1}},$$

then it follows from (10.7) that, uniformly for $\alpha \in I_c'$ and $z \in B_2$ with $|z| = \beta(\alpha)$,

$$(10.10) \qquad G_{n,\alpha}(z) \lesssim |H_n(z)| (\log n)^{1/2}.$$

Let $z \in B_2$ and let $\delta_n = 1/\log n$. Since $W_n$ is analytic,

$$(10.11) \qquad |W_n(z)|^2 \leq \frac{1}{\pi \delta_n^2} \int_{|w-z| < \delta_n} |W_n(z)|^2 \, dm(w),$$

where $m$ is the two-dimensional Lebesgue measure.

Further, if $|w - z| < \delta_n$, then $n^{\lambda_1(|z|)-1}$ and $n^{\lambda_1(|w|)-1}$ differ by at most a constant factor. Hence,

$$(10.12) \qquad |H_n(z)|^2 \lesssim (\log n)^2 \int_{|w-z| < \delta_n} |H_n(z)|^2 \, dm(w).$$

Fix $\delta > 0$ (we will be more precise below) and let $B_2^\delta = \{w : \mathrm{dist}(w, B_2) \leq \delta\}$. Then, by (10.12), if $n$ is sufficiently large that $\delta_n < \delta$, we have

$$(10.13) \qquad \sup_{z \in B_2} |H_n(z)|^2 \lesssim (\log n)^2 \int_{B_2^\delta} |H_n(z)|^2 \, dm(z).$$

We now use Lemma 8.4. By Lemma 8.1 and the assumption $I_c \subset I$,

$$\max(\lambda_1(|z|^2) - 1, 2\Re(\lambda_1(z)) - 2, 0) < 2\lambda_1(|z|) - 2$$

on $B_2$. By continuity, this also holds on $B_2^\delta$ if $\delta > 0$ is small enough and, by Lemma 8.4 and compactness, there exists $\eta > 0$ such that

$$\mathbb{E}|W_n(z)|^2 \lesssim n^{2\lambda_1(|z|)-2-\eta}$$

uniformly for $z \in B_2^\delta$. Thus by (10.9),

$$\mathbb{E}|H_n(z)|^2 \lesssim n^{-\eta}$$

and by (10.13),

$$\mathbb{E}\Big(\sup_{z \in B_2} |H_n(z)|^2\Big) \lesssim (\log n)^2 \int_{B_2^\delta} \mathbb{E}|H_n(z)|^2 \, dm(z) \lesssim (\log n)^2 n^{-\eta}.$$

Hence, $\sup_{z \in B_2} \log n |H_n(z)| \xrightarrow{\text{P}} 0$ and the result follows by (10.10). $\square$

PROOF OF THEOREM 1.1. By Theorem 1.2, $W_n(z)/\mathbb{E}W_n(z) \xrightarrow{\text{d}} Y$ in $\mathcal{H}(B)$ and thus for every compact subset $B_c$ in the space $C(B_c)$ with the uniform topology. Hence, Lemma 10.1(ii) and [2], Theorem 5.5, imply that

$$T_n(W_n/\mathbb{E}W_n) \xrightarrow{\text{d}} Y \qquad \text{in } D(I_c').$$

Finally, Lemma 10.2 and (10.8) imply that, provided $I_c \subset I$,

$$\frac{X_{n,\lfloor \alpha \log n \rfloor}}{\mathbb{E}X_{n,\lfloor \alpha \log n \rfloor}} - T_n\Big(\frac{W_n(z)}{\mathbb{E}W_n(z)}\Big) \xrightarrow{\text{P}} 0,$$

uniformly on $I_c'$. The theorem then follows. $\square$

REMARK 10.1. The reason that we have to restrict ourselves to the interval $I'$ (and cannot extend our result to a larger interval; compare with the discussion of the critical values in Section 1) is that we use an $L^2$-estimate for $|W_n(z)|$ in the proof of Lemma 10.2 that only works if $\alpha \in I'$. In fact, $I'$ is the largest interval in which we have $L_2$-convergence to the process $Y(z)$.

However, it is very likely that one can prove similar estimates for $\mathbb{E}|W_n(z)|^p$ for any $p > 1$ and that our method of proof, using the version of (10.11) for $p$th powers, would then prove Theorem 1.1 for the largest possible interval $(\alpha_0, \alpha_+)$.

## 11. The external profile.
In this final section, we will discuss a variation of Theorem 1.1 which concerns a similarly defined profile process.

The *external profile* $Y_{n,k}$ denotes the number of *free positions* at level $k$ in a tree with $n$ keys. A free position is a position where the $(n+1)$st key can be placed, for example, $Y_{0,0} = m - 1$ and $Y_{0,k} = 0$ for $k \geq 1$. More precisely, we have

$$(11.1) \qquad Y_{n,k} \overset{\text{d}}{=} Y_{V_{n,1},k-1}^{(1)} + Y_{V_{n,2},k-1}^{(2)} + \cdots + Y_{V_{n,m},k-1}^{(m)},$$

jointly in $k \geq 0$ for every $n \geq m - 1$, where the random vector $\mathbf{V}_n = (V_{n,1}, V_{n,2}, \ldots, V_{n,m})$ is as in Section 2 and is the same for every $k \geq 0$, and $\mathbf{Y}_n^{(j)} = (Y_{n,k}^{(j)})_{k \geq 0}$, $j = 1, \ldots, m$, are independent copies of $\mathbf{Y}_n$ that are also independent of $\mathbf{V}_n$. The initial conditions are $Y_{n,0} = 0$ for $n \geq m - 1$, and for $n \leq m - 2$, we simply have $Y_{n,0} = m - 1 - n$ for $0 \leq n \leq m - 2$ and $Y_{n,k} = 0$ for $k \geq 1$. However, we should mention that the initial conditions for $Y_{0,0}, Y_{1,0}, \ldots, Y_{m-2,0}$ only affect implicit constants in our analysis. The limit theorem (Theorem 11.1) is not affected. For example, one can also use $Y_{n,0} = 1$ for $0 \leq n \leq m-2$ in order to count the number of nodes where one can place a new item.

Let $U_n(z) = \sum_k Y_{n,k} z^k$ denote the random *external profile polynomial*. By (11.1), it is recursively given by $U_n(z) = m - 1 - n$ for $0 \leq n \leq m-2$ and

$$(11.2) \qquad U_n(z) \overset{\mathrm{d}}{=} z U_{V_{n,1}}^{(1)}(z) + z U_{V_{n,2}}^{(2)}(z) + \cdots + z U_{V_{n,m}}^{(m)}(z), \qquad n \geq m-1,$$

where $U_\ell^{(j)}(z)$, $j = 1, \ldots, m$, are independent copies of $U_\ell(z)$ that are independent of $\mathbf{V}_n$, $\ell \geq 0$. From this relation, we obtain (similarly to the above) a recurrence for the expected external profile polynomial $\mathbb{E}U_n(z)$. We have, using (2.2), for $n \geq mt + m - 1$,

$$(11.3) \qquad \mathbb{E}U_n(z) = mz \sum_{\ell=0}^{n-1} \frac{\binom{\ell}{t}\binom{n-\ell-1}{(m-1)t+m-2}}{\binom{n}{mt+m-1}} \mathbb{E}U_\ell(z).$$

By [9, Theorem 1(i)], we obtain, as above,

$$(11.4) \qquad \mathbb{E}U_n(z) \sim \overline{E}(z) n^{\lambda_1(z)-1}$$

for some analytic function $\overline{E}(z)$ with $\overline{E}(z) > 0$ for $z > 0$. Moreover, this limit relation is true for all $z \in D_{-t}$, not only for $z \in D_1$, as we will see in a moment. Since $D_1 \subset D_{-t}$, we can expect that corresponding limit theorems hold for a larger range.

The fact that (11.4) holds for $z \in D_{-t}$ needs some explanation. If we just use [9], Theorem 1(i), this creates the impression that $z \in D_0$ is the largest region for (11.4) since [9], Theorem 1(i) assumes that $\Re(\lambda_1(z)) > 0$. Furthermore, the indicial polynomial $\Lambda(\lambda; z)$ always has the roots $0, -1, -2, \ldots, -t+1$, that is, if $\Re(\lambda_1(z)) \leq 0$, then the dominant root of $\Lambda(\lambda; z)$ is always 0.

However, the contribution of a simple root $\lambda(z)$ of $\Lambda(\lambda; z)$ to the behavior of $\mathbb{E}U_n(z)$ is of the form $c(z)(-1)^n \binom{-\lambda(z)}{n}$. This implies that the roots $0, -1, -2, \ldots, -t+1$, provided they are simple, only matter for $n < t$, that is, they do not affect the asymptotics of $\mathbb{E}U_n(z)$. Hence, if $\lambda_1(z)$, the dominant root of $F(\theta) = z$, is different from $0, -1, -2, \ldots, -t+1$, then (11.4) is also true. Further, by Theorem A.1, (11.4) holds uniformly for any compact set contained in $D_{-t}$ with $\lambda_1(z) \notin \{0, -1, -2, \ldots, -t+1\}$.

Nevertheless, a little bit of more careful analysis reveals that the exceptional values $0, -1, -2, \ldots, -t+1$ are only present in the analysis, but not in the asymptotic result. The limit relation (11.4) *extends continuously* to all $z \in D_{-t}$. For example, if $t \geq 1$ and $\lambda_1(z)$ is close to zero, then it follows that the power series of $\mathbb{E}U_n(z)$ can be represented as (compare to the Appendix)

$$(11.5) \qquad \sum_{n \geq 0} \mathbb{E}U_n(z)\zeta^n = \Gamma(\lambda_1(z))\overline{E}(z)\left(\frac{1}{(1-\zeta)^{\lambda_1(z)}} - 1\right)$$

$$+ \text{ smaller order terms.}$$

Further, if $\lambda_1(z_0) = 0$, then $\Lambda(\lambda; z)$ has a double zero and we have

$$(11.6) \qquad \sum_{n \geq 0} \mathbb{E}U_n(z_0)\zeta^n = \overline{E}(z_0) \log \frac{1}{1-\zeta} + \text{smaller order terms.}$$

Note that (11.5) and (11.6) are *consistent* for $z \to z_0$ and they imply that we also have (11.4) uniformly in a neighborhood of $z = z_0$. Similar phenomena appear if $\lambda_1(z)$ is close to $-1, -2, \ldots, -t+1$.

It can also be shown, using (an analytic extension of) the formula in [9], Theorem 1(i), as for $E(z)$ in Lemma 3.1, that $\overline{E}(z) > 0$ if $z > 0$ is real.

The expected external profile $\mathbb{E}Y_{n,k}$ is, thus, given by the following asymptotic formula which can be proven exactly as Lemma 8.3.

LEMMA 11.1. *Suppose that $\alpha_1, \alpha_2$ with $0 < \alpha_1 < \alpha_2 < \infty$ are given and let $\beta(\alpha)$ be defined by $\beta(\alpha)\lambda_1'(\beta(\alpha)) = \alpha$. Then*

$$\mathbb{E}Y_{n,k} = \frac{\overline{E}(\beta(\alpha_{n,k}))n^{\lambda_1(\beta(\alpha_{n,k})) - \alpha_{n,k}\log(\beta(\alpha_{n,k})) - 1}}{\sqrt{2\pi(\alpha_{n,k} + \beta(\alpha_{n,k})^2\lambda_1''(\beta(\alpha_{n,k})))\log n}}(1 + O((\log n)^{-1/2}))$$

*uniformly for $\alpha_{n,k} = k/\log n \in [\alpha_1, \alpha_2]$ as $n, k \to \infty$.*

Hence, we can proceed as above and obtain the following variation of Theorem 1.1, with some minor differences in the proofs. [For example, we use the fact that $\mathbb{E}(1 + V_{n,r})^a = O(n^a)$ for any $a > -t - 1$; we omit the verification of this estimate.]

THEOREM 11.1. *Let $m \geq 2$ and $t \geq 0$ be given integers and let $(Y_{n,k})_{k\geq 0}$ be the external profile of the corresponding random search tree with $n$ keys.*

*Let $J = \{\beta > 0 : \lambda_1(\beta^2) < 2\lambda_1(\beta) - 1\}$ and $J' = \{\beta\lambda_1'(\beta) : \beta \in J\}$. We then have*

$$(11.7) \qquad \left(\frac{Y_{n,\lfloor\alpha\log n\rfloor}}{\mathbb{E}Y_{n,\lfloor\alpha\log n\rfloor}}, \alpha \in J'\right) \xrightarrow{\mathrm{d}} (Y(\beta(\alpha)), \alpha \in J'),$$

*where $Y(z)$ is as in Theorem 1.1.*

The difference between Theorems 1.1 and 11.1 is that Theorem 11.1 is true for a larger range for $k/\log n$ since $I' \subset J'$. The reason is that the internal profile $X_{n,k}$ has a phase transition at level $\alpha_0$ that is not present for the external profile (compare with the discussion of critical constants in the Introduction).

Further, note that we can also deal with the process $R_{n,k} := (m-1)m^k - X_{n,k}$ that can be approximated by $Y(z)$ in the range $\alpha \in (0, \alpha_0) \cap J'$. (We do not work out the details here.)

We close this section with a vector-valued generalization of Theorem 1.1. As we know from the introduction, a node in an $m$-ary search tree stores one or several of the keys up to at most $m-1$. This means that we can partition the nodes into types. We say that a node has *type $j$* ($j \in \{1, 2, \ldots, m-1\}$) if it stores exactly $j$ keys. Further, we can also extend this definition to $j = 0$ if we define nodes of type $0$ as external nodes.

Now, let $X_{n,k}^{(j)}$ denote the number of nodes of type $j$ at level $k$ in a random $m$-ary search tree with $n$ keys. We can then prove the following theorem.

THEOREM 11.2. *Let $m \geq 2$ and $t \geq 0$ be given integers and let $(X_{n,k}^{(0)}, \ldots, X_{n,k}^{(m-1)})$, $k \geq 0$, be the random profile vector of the random search tree with $n$ keys.*

*Let $I = \{\beta > 0 : 1 < \lambda_1(\beta^2) < 2\lambda_1(\beta) - 1\}$ and $I' = \{\beta\lambda_1(\beta) : \beta \in I\}$. We then have*

$$(11.8) \qquad \begin{aligned} &\left( \frac{X_{n,\lfloor \alpha \log n \rfloor}^{(0)}}{\mathbb{E} X_{n,\lfloor \alpha \log n \rfloor}^{(0)}}, \ldots, \frac{X_{n,\lfloor \alpha \log n \rfloor}^{(m-1)}}{\mathbb{E} X_{n,\lfloor \alpha \log n \rfloor}^{(m-1)}} ; \alpha \in I' \right) \\ &\qquad \overset{\mathrm{d}}{\longrightarrow} (Y(\beta(\alpha)), \ldots, Y(\beta(\alpha)) ; \alpha \in I'). \end{aligned}$$

We do not work out the details, but the same proof techniques as for the proof of Theorem 1.1 also work here, using Remark 7.1. We obtain convergence to a vector $(Y_i(z))_{i=0}^{m-1}$ satisfying

$$(Y_i(z))_i \overset{\mathrm{d}}{=} (z V_1^{\lambda_1(z)-1} Y_i^{(1)}(z) + z V_2^{\lambda_1(z)-1} Y_i^{(2)}(z) + \cdots + z V_m^{\lambda_1(z)-1} Y_i^{(m)}(z))_i,$$

where $(Y_i^{(j)}(z))_i$, $1 \leq j \leq m$, are independent copies of $(Y_i(z))_i$, and this equation is solved by $Y_0(z) = \cdots = Y_{m-1}(z) = Y(z)$. It follows that the profiles $X_{n,k}^{(i)}$, $0 \leq i \leq m-1$, are asymptotically proportional.

Note that (11.8) is a functional limit result. Hence, joint convergence for several different arguments $\alpha$ follows. In particular, we obtain limits for vectors $(X_{n,\lfloor \alpha_0 \log n \rfloor}^{(0)}, \ldots, X_{n,\lfloor \alpha_{m-1} \log n \rfloor}^{(m-1)})$ with fixed $\alpha_0, \ldots, \alpha_{m-1}$.

Note that Theorem 1.1 follows directly from Theorem 11.2, since

$$X_{n,k} = X_{n,k}^{(1)} + 2X_{n,k}^{(2)} + \cdots + (m-1)X_{n,k}^{(m-1)}.$$

# APPENDIX A: DETAILED PROOFS OF LEMMAS 8.2 AND 8.4

In Lemmas 8.2 and 8.4, we used results by Chern, Hwang and Tsai [9], where we claimed uniformity in $z$ (in certain compact sets). This uniformity can be verified by a tedious checking of the proofs in [9] (Hwang, personal communication), but, for completeness, we give a detailed proof here; see also Chern and Hwang [8] for the case $m = 2$.

Fortunately, we are in the situation where the relevant generating functions can be analytically continued outside the unit disc so that the singularity analysis of Flajolet and Odlyzko [13] (see also Flajolet and Sedgewick [14], Chapter 6) applies. As pointed out in [9], page 197, this simplifies the arguments considerably, so we consider only this case.

We introduce the generating function $\Psi(\zeta; z) := \sum_{n \geq 0} \mathbb{E} W_n(z) \zeta^n$. Let $\Lambda(\theta; z)$ be the polynomial (in $\theta$) of degree $r := mt + m - 1$ defined in (3.4) and let $\vartheta$ denote the differential operator $(1 - \zeta) \frac{d}{d\zeta}$. Further, let $b_n(z) := \mathbb{E} W_n(z)$ for $n < r$ and $b_n(z) := m - 1$ for $n \geq r$, and define the generating function $g(\zeta; z) := \sum_{n \geq 0} b_n(z) \zeta^n$. Then, as is shown in [9], (3.3) is equivalent to the differential equation (in $\zeta$, with $z$ fixed)

$$(A.1) \qquad \Lambda(\vartheta; z) \Psi(\zeta; z) = \phi(\zeta; z) := (1 - \zeta)^r \frac{\partial^r}{\partial \zeta^r} g(\zeta; z).$$

Note that $g(\zeta; z)$, for every $z$, differs from $(m - 1)(1 - \zeta)^{-1}$ by a polynomial in $\zeta$ of degree less than $r$. Hence,

$$(A.2) \quad \phi(\zeta; z) := (m - 1)(1 - \zeta)^r \frac{\partial^r}{\partial \zeta^r} (1 - \zeta)^{-1} = (m - 1)r!(1 - \zeta)^{-1}.$$

We study solutions to (A.1) in some generality and state a theorem where $\Lambda$ and $\phi$ can be rather arbitrary polynomials and analytic functions, respectively, depending on a parameter $z$. (The parameter set $K$ can be any set, although we only need subsets of the complex plane for the present paper.)

A $\Delta$-*domain* is a domain of the type

$$\Delta(R, \delta) := \{z \in \mathbb{C} : |z| < R \text{ and } |\arg(z - 1)| > \pi/2 - \delta\}$$

for some $R > 1$ and $\delta \in (0, \pi/2)$.

THEOREM A.1. *Let* $r \geq 1$ *and for each* $z \in K$, *let* $\Lambda_z(\theta)$ *be a monic polynomial in* $\theta$ *of degree* $r$, *with coefficients that are bounded functions of* $z$. *Moreover, let* $\phi_z(\zeta)$ *be an analytic function of* $\zeta$ *in the unit disc for every* $z \in K$ *and let* $\Psi_z(\zeta)$ *be a formal power series that solves the differential equation*

$$(A.3) \qquad \Lambda_z(\vartheta) \Psi_z(\zeta) = \phi_z(\zeta).$$

*We denote the roots of* $\Lambda_z(\lambda) = 0$ *(counted with multiplicities) by* $\lambda_j(z)$, $j = 1, \ldots, r$, *arranged in decreasing order of their real parts:* $\Re \lambda_1(z) \geq \Re \lambda_2(z) \geq \cdots$.

(i) *Assume that, for each* $z \in K$, $\phi_z(\zeta)$ *extends to an analytic function in a fixed* $\Delta$-*domain* $\Delta = \Delta(R, \delta)$ *and that for some constants* $\alpha \in (-\infty, \infty)$, $d \geq 0$ *and* $\eta \in (0, 1/2)$, *uniformly in all* $z \in K$ *and* $\zeta \in \Delta$,

$$(A.4) \qquad \phi_z(\zeta) = \begin{cases} O(1), & |1 - \zeta| \geq \eta, \\ O(|1 - \zeta|^{-\alpha} |\log|1 - \zeta||^d), & |1 - \zeta| \leq \eta, \end{cases}$$

*and that, again uniformly in $z \in K$,*

$$\text{(A.5)} \qquad \frac{\partial^k}{\partial \zeta^k} \Psi_z(0) = O(1), \qquad k = 0, \ldots, r-1.$$

*Each $\Psi_z(\zeta)$ then converges in the unit disc and extends to an analytic function in $\Delta$ such that, uniformly in all $z \in K$ and $\zeta \in \Delta$,*

$$\Psi_z(\zeta) = \begin{cases} O(1), & |1-\zeta| \geq \eta, \\ O(|1-\zeta|^{-(\alpha \vee \Re \lambda_1(z))} |\log|1-\zeta||^{d+r}), & |1-\zeta| \leq \eta. \end{cases}$$

(ii) *Assume, further, that $\varepsilon > \varepsilon_1 > 0$ and that, for all $z \in K$, $\Re \lambda_2(z) \leq \Re \lambda_1(z) - \varepsilon$ and $\alpha \leq \Re \lambda_1(z) - \varepsilon$. (The first assumption is void if $r = 1$.) Then, for some function $c(z)$, uniformly in all $z \in K$ and $\zeta \in \Delta$,*

$$\Psi_z(\zeta) = c(z)(1-\zeta)^{-\lambda_1(z)} + O(|1-\zeta|^{\varepsilon_1 - \Re \lambda_1(z)}).$$

*In particular, this holds for some $\varepsilon_1 > 0$ if $K$ is a compact topological space, the coefficients of $\Lambda_z$ are continuous functions of $z$ and $\Re \lambda_1(z) > \Re \lambda_2(z)$ and $\Re \lambda_1(z) > \alpha$ for each $z \in K$.*

PROOF. We have $\Lambda_z(\vartheta) = \prod_{i=1}^r (\vartheta - \lambda_1(z))$. The roots $\lambda_i(z)$ are not always continuous functions of the coefficients of $\Lambda_z$ (because of ambiguity in labeling the roots), but $\max_i |\lambda_i(z)|$ is, and since the coefficients are bounded, it follows that $\sup\{|\lambda_i(z)| : z \in K, i = 1, \ldots, r\} < \infty$. We may thus treat the factors $\vartheta - \lambda_i(z)$ one by one and, by induction, it suffices to prove part (i) for the case $r = 1$, $\Lambda_z(\vartheta) = \vartheta - \lambda(z)$, where $\lambda(z)$ is bounded, provided we show that the bounds are also uniform in $\alpha \in A$ for any bounded set $A$.

In this case (see [9], Lemma 1), it is easily seen that, for each $y$, there is a unique power series $\Psi_z$ satisfying (A.3) with $\Psi_z(0) = y$. Moreover, the solution $\Psi_z$ is given by the analytic function

$$\text{(A.6)} \quad \Psi_z(\zeta) = \Psi_z(0)(1-\zeta)^{-\lambda(z)} + (1-\zeta)^{-\lambda(z)} \int_0^\zeta (1-w)^{\lambda(z)-1} \phi_z(w) \, dw.$$

Clearly, (A.6) defines $\Psi_z$ as an analytic function in $\Delta$ because $\phi_z$ is such a function.

For (i), it only remains to estimate this solution for $\zeta \in \Delta$, with $\alpha = O(1)$ and $\lambda = \lambda(z) = O(1)$. We first note that the value of $\eta$ is immaterial, so we may, for convenience, assume that $\eta < R - 1$. For $|1-\zeta| \geq \eta$, we can choose an integration path in (A.6) of bounded length and contained in the region $|1-w| \geq \eta$, and it follows that $\Psi_z(\zeta) = O(1)$.

Now, assume that $\zeta \in \Delta$ and $|1-\zeta| \leq \eta$. Let $\gamma := \arg(1-\zeta)$ and $\zeta' := 1 - \eta e^{i\gamma}$. Then

$$\int_0^\zeta (1-w)^{\lambda-1} \phi_z(w) \, dw = \int_0^{\zeta'} (1-w)^{\lambda-1} \phi_z(w) \, dw + \int_{\zeta'}^\zeta (1-w)^{\lambda-1} \phi_z(w) \, dw$$

$$\lesssim 1 + \int_{|1-\zeta|}^{\eta} x^{\Re\lambda-1} x^{-\alpha} |\log x|^d \, dx$$

$$\lesssim |\log|1-\zeta||^d \int_{|1-\zeta|}^{1} x^{-(\alpha-\Re\lambda)-1} \, dx.$$

For any real $\beta \neq 0$ and $0 < y < 1$, by the mean value property for $t \mapsto e^{t \log(1/y)}$,

$$\int_y^1 x^{-\beta-1} \, dx = \frac{1-y^{-\beta}}{-\beta} = \frac{e^{\beta \log(1/y)} - 1}{\beta}$$

$$\leq |\log y| e^{(\beta \vee 0) \log(1/y)} = |\log y| y^{-(\beta \vee 0)}.$$

This is evidently also true for $\beta = 0$ and thus

$$\int_0^{\zeta} (1-w)^{\lambda-1} \phi_z(w) \, dw \lesssim |\log|1-\zeta||^{d+1} |1-\zeta|^{-((\alpha-\Re\lambda)\vee 0)}.$$

Consequently, (A.6) yields, recalling $|1-\zeta| \leq \eta \leq 1/2$,

$$|\Psi_z(\zeta)| \lesssim (1+|\log|1-\zeta||^{d+1} |1-\zeta|^{-((\alpha-\Re\lambda)\vee 0)}) |(1-\zeta)^{-\lambda}|$$

$$\lesssim |\log|1-\zeta||^{d+1} |1-\zeta|^{-(\alpha\vee\Re\lambda)}.$$

For part (ii), we factorize $\Lambda_z(\vartheta) = \Lambda_z^*(\vartheta)(\vartheta - \lambda_1(z))$ and let $\Psi_z^*(\zeta) = (\vartheta - \lambda_1(z))\Psi_z(\zeta)$. Thus, $\Lambda_z^*(\vartheta)\Psi_z^*(\zeta) = \phi_z(\zeta)$. By (i) applied to $\Lambda_z^*$ (or directly if $r=1$), for $|1-\zeta| \leq \eta$,

$$(A.7) \qquad \Psi_z^*(\zeta) = O(|1-\zeta|^{-(\alpha\vee\Re\lambda_2(z))} |\log|1-\zeta||^{d+r-1})$$

$$= O(|1-\zeta|^{-\Re\lambda_1(z)+\varepsilon_1}).$$

We now use (A.6), with $\lambda(z) = \lambda_1(z)$ and $\phi_z$ replaced by $\Psi_z^*$. By (A.7), the integral $c_1(z) := \int_0^1 (1-w)^{\lambda_1(z)-1} \Psi_z^*(w) \, dw$ converges and if we define $c(z) := \Psi_z(0) + c_1(z)$, then, again using (A.7) we have

$$\Psi_z(\zeta) - c(z)(1-\zeta)^{-\lambda_1(z)}$$

$$= -(1-\zeta)^{-\lambda_1(z)} \int_{\zeta}^1 (1-w)^{\lambda_1(z)-1} \Psi_z^*(w) \, dw$$

$$\lesssim |1-\zeta|^{-\Re\lambda_1(z)} \int_0^{|1-\zeta|} x^{\Re\lambda_1(z)-1-\Re\lambda_1(z)+\varepsilon_1} \, dx$$

$$\lesssim |1-\zeta|^{-\Re\lambda_1(z)+\varepsilon_1},$$

provided $|1-\zeta| \leq \eta$. The case $|1-\zeta| \geq \eta$ follows immediately from $\Psi_z(\zeta) = O(1)$.

The final statement follows because, under these assumptions, $\Re\lambda_1(z)$ and $\Re\lambda_2(z)$ are continuous functions of $z \in K$, so the compactness of $K$ yields the existence of an $\varepsilon > 0$ such that $\Re\lambda_1(z) - \Re\lambda_2(z) \geq \varepsilon$ and $\Re\lambda_1(z) - \alpha \geq \varepsilon$ for all $z \in K$. $\square$

PROOF OF LEMMA 8.2. We use Theorem A.1 with $\Lambda_z(\zeta) = \Lambda(\zeta; z)$, $\Psi_z(\zeta) = \Psi(\zeta; z)$ and $\phi_z(\zeta) = \phi(\zeta; z)$, as above. [By (A.2), $\phi_z(\zeta)$ is actually independent of $z$.] We see from (A.2) that (A.4) holds with $\alpha = 1$ and $d = 0$. Hence, Theorem A.1(i) yields, for $z \in K$, $\zeta \in \Delta$,

$$\Psi(\zeta; z) = \begin{cases} O(1), & |1 - \zeta| \geq \eta, \\ O(|1 - \zeta|^{-(\Re\lambda_1(z)\vee 1)}|\log|1 - \zeta||^r), & |1 - \zeta| \leq \eta, \end{cases}$$

and Lemma 8.2(ii) follows by standard singularity analysis (see, e.g., [14], Chapter 6).

Similarly, Lemma 8.2(i) follows easily from Theorem A.1(ii). The constant $E(z)$ is necessarily the same as in Lemma 3.1. $\square$

PROOF OF LEMMA 8.4. This time, we use (8.12). Let $a_n(z) := \mathbb{E}|W_n(z)|^2$ and let $b_n(z)$ be the $O$ term in (8.12). (8.12) then implies that the analogue of (A.1) holds for the generating functions of $a_n$ and $b_n$, with $\Lambda(\zeta, z)$ replaced by $\Lambda(\zeta, |z|^2)$. However, it is not clear that these generating functions extend to a $\Delta$-domain. Therefore, we instead take $g_z(\zeta) := Ch(\zeta)$, where

$$h(\zeta) := (1 - \zeta)^{-\max\{2\Re(\lambda_1(z)) - 1, 1\}}(-\log(1 - \zeta)/\zeta)^{2D}$$

and $C$ is a constant chosen sufficiently large that $|b_n(z)| \leq [\zeta^n]g(\zeta)$ for all $z \in K$ and $n \geq 0$. This is possible because $[\zeta^n]h(\zeta) = \Theta(n^{\max\{2\Re(\lambda_1(z)) - 2, 0\}}(\log n)^{2D})$. We then let $\Lambda_z(\zeta) := \Lambda(\zeta, |z|^2)$ and define $\Psi_z(\zeta)$ to be the solution of (A.3) with initial conditions $[\zeta^n]\Psi_z(\zeta) = \mathbb{E}|W_n(z)|^2$ for $n < r$. It now follows from (8.12) and induction that $\mathbb{E}|W_n(z)|^2 \leq [\zeta^n]\Psi_z(\zeta)$ for all $n$, so it suffices to estimate $[\zeta^n]\Psi_z(\zeta)$.

Clearly, (A.4) is satisfied with $\alpha = \max\{2\Re(\lambda_1(z)) - 1, 1\}$ and $d = 2D$. Theorem A.1(i) thus yields, uniformly for $z \in K$,

$$\Psi_z(\zeta) = \begin{cases} O(1), & |1 - \zeta| \geq \eta, \\ O(|1 - \zeta|^{-\max\{\lambda_1(|z|^2), 2\Re(\lambda_1(z)) - 1, 1\}}|\log|1 - \zeta||^{2D+r}), \\ & |1 - \zeta| \leq \eta. \end{cases}$$

Standard singularity analysis then yields

$$[\zeta^n]\Psi_z(\zeta) = O(n^{\max\{\lambda_1(|z|^2) - 1, 2\Re(\lambda_1(z)) - 2, 0\}}(\log n)^{2D+r})$$

and (8.10) follows. $\square$

**Acknowledgment.** We thank Ludger Rüschendorf for helpful comments.

M. DRMOTA
INSTITUTE OF DISCRETE MATHEMATICS
  AND GEOMETRY
TU WIEN
WIEDNER HAUPTSTR. 8-10
A-1040 WIEN
AUSTRIA
E-MAIL: michael.drmota@tuwien.ac.at

S. JANSON
DEPARTMENT OF MATHEMATICS
UPPSALA UNIVERSITY
PO BOX 480
SE-751 06 UPPSALA
SWEDEN
E-MAIL: svante.janson@math.uu.se

R. NEININGER
DEPARTMENT OF MATHEMATICS
  AND COMPUTER SCIENCE
J. W. GOETHE UNIVERSITY
D-60325 FRANKFURT A.M.
GERMANY
E-MAIL: neiningr@math.uni-frankfurt.de